\documentclass[a4paper,notitlepage, eqno]{article}%
\usepackage{amssymb}
\usepackage{graphicx}
\usepackage{amsmath}
\usepackage{amsthm}
\usepackage{amsfonts}%
\providecommand{\U}[1]{\protect\rule{.1in}{.1in}}
\theoremstyle{plain}
\newtheorem{theorem}{Theorem}[section]

\newtheorem{cor}[theorem]{Corollary}

\newtheorem{lem}[theorem]{Lemma}

\newtheorem{prop}[theorem]{Proposition}

\theoremstyle{definition}

\newtheorem{rem}[theorem]{Remark}
\numberwithin{equation}{section}
\numberwithin{theorem}{section}

\ifx\pdfoutput\relax\let\pdfoutput=\undefined\fi
\newcount\msipdfoutput
\ifx\pdfoutput\undefined\else
\ifcase\pdfoutput\else
\ifx\paperwidth\undefined\else
\ifdim\paperheight=0pt\relax\else\pdfpageheight\paperheight\fi
\ifdim\paperwidth=0pt\relax\else\pdfpagewidth\paperwidth\fi
\fi\fi\fi
\begin{document}

\title{Nonnegative solutions of an indefinite sublinear Robin problem I: positivity,
exact multiplicity, and existence of a subcontinuum. 
 }
\author{Uriel Kaufmann\thanks{FaMAF, Universidad Nacional de C\'{o}rdoba, (5000)
C\'{o}rdoba, Argentina. \textit{E-mail address: }kaufmann@mate.uncor.edu. Partially supported by Secyt-UNC 33620180100016CB.} ,
Humberto Ramos Quoirin \thanks{Departamento de Matem\'aticas y C.C., Universidad de Santiago de Chile, Casilla 307,
Correo 2, Santiago, Chile. \textit{E-mail address: }humberto.ramos@usach.cl. Supported by Fondecyt grants 1161635, 1171532, 1171691, 1181825}  \thanks{CIEM-FaMAF, Universidad Nacional de C\'{o}rdoba, (5000)
C\'{o}rdoba, Argentina.} , 
Kenichiro Umezu\thanks{Department of Mathematics, Faculty of Education,
Ibaraki University, Mito 310-8512, Japan. \textit{E-mail address:
}kenichiro.umezu.math@vc.ibaraki.ac.jp. Supported by JSPS KAKENHI Grant Numbers 15K04945 and 18K03353.}
\and \noindent}
\maketitle

\begin{abstract}
Let $\Omega\subset\mathbb{R}^{N}$ ($N\geq1$) be a smooth bounded domain,
$a\in C(\overline{\Omega})$ a sign-changing function, and $0\leq q<1$. We investigate the
Robin problem
\[
\begin{cases}
-\Delta u=a(x)u^{q} & \mbox{in $\Omega$},\\
u\geq0 & \mbox{in $\Omega$},\\
\partial_{\nu}u=\alpha u & \mbox{on $\partial \Omega$},
\end{cases}
\leqno{(P_{\alpha})}
\]
where  $\alpha\in\lbrack-\infty,\infty)$ and $\nu$ is the unit
outward normal to $\partial\Omega$. Due to the lack of strong maximum
principle structure, this problem may have \textit{dead
core} solutions. However, for a large class of weights $a$ we recover a
\textit{positivity} property when $q$ is close to $1$, which considerably
simplifies the structure of the solution set.  Such property, combined with a
bifurcation analysis and a suitable change of variables, enables us to show
the following exactness result for these values of $q$: $(P_{\alpha})$ has
\textit{exactly} one nontrivial solution for $\alpha\leq0$, \textit{exactly}
two nontrivial solutions for $\alpha>0$ small, and \textit{no} such solution
for $\alpha>0$ large. Assuming some further conditions on $a$, we show that
these solutions lie on a subcontinuum. These results rely partially on (and
extend) our previous work \cite{KRQU16}, where the cases $\alpha=-\infty$
(Dirichlet) and $\alpha=0$ (Neumann) have been considered. We also obtain some
results for arbitrary $q\in\left[  0,1\right)  $. Our approach combines
mainly bifurcation techniques, the sub-supersolutions method, and \textit{a
priori} lower and upper bounds.

\end{abstract}

%----------------

%\author{Uriel Kaufmann\thanks{FaMAF, Universidad Nacional de C\'{o}rdoba, (5000) C\'{o}rdoba, Argentina. \textit{E-mail address: }kaufmann@mate.uncor.edu} , Humberto Ramos Quoirin\thanks{Universidad de Santiago de Chile, Casilla 307, Correo 2, Santiago, Chile. \textit{E-mail address: }humberto.ramos@usach.cl}, Kenichiro Umezu\thanks{Department of Mathematics, Faculty of Education, Ibaraki University, Mito 310-8512, Japan. \textit{E-mail address: }kenichiro.umezu.math@vc.ibaraki.ac.jp} \and \noindent}

%\begin{abstract}

%\end{abstract}

%restart on 20 feb 2018.

%\begin{description}
%\item[modification] \textit{Following the manner in ANS paper, the term \textbf{subcontinuum} is used for a closed and connected subset, whereas the term \textbf{component} is used for a maximal closed, connected subset. }
%\end{description}

%----------------------------------------------

\section{Introduction}

\label{sec:Intro}

This article is devoted to a class of indefinite elliptic pdes, whose
prototype is the equation
\[
-\Delta u=a(x)u^{q}\quad\text{in} \quad\Omega,
\]
where $\Omega\subset\mathbb{R}^{N}$ ($N\geq1$) is a bounded and smooth domain,
and $a\in C(\overline{\Omega})$ is a sign-changing function. Over the past  
decades, many works have addressed basic issues on nonnegative solutions of this equation (under different boundary conditions) in the superlinear
case $q>1$ \cite{AT,ALG,BCN,BGH,LG,Ou,Te}. On the other hand, much less
attention has been given to the sublinear problem, i.e. with $0<q<1$, which
will be considered here. In particular, we shall highlight the main contrasts
between these two cases.

We consider nonnegative solutions of the above equation under a Robin boundary
condition, i.e. the problem:
\[%
\begin{cases}
-\Delta u=a(x)u^{q} & \mbox{in $\Omega$},\\
u\geq0 & \mbox{in $\Omega$},\\
\partial_{\nu}u=\alpha u & \mbox{on $\partial \Omega$}.
\end{cases}
\leqno{(P_{\alpha})}
\]
Here $\nu$ is the unit outward normal to $\partial\Omega$, $\partial_{\nu
}:=\frac{\partial}{\partial\nu}$, and $\alpha\in\lbrack-\infty,\infty)$. When
$\alpha=-\infty$ the boundary condition is understood as $u=0$ on
$\partial\Omega$, so that we treat in particular the Dirichlet ($\alpha
=-\infty$) and Neumann ($\alpha=0$) problems.

Our main interest is the structure of the solutions set of this problem. By a
\textit{solution} of
%\marginpar{modified}
$(P_{\alpha})$ we mean a nonnegative function $u\in W^{2,r}\left(
\Omega\right)  $, with $r>N$,
%\marginpar{added}
that satisfies the equation for the weak derivatives and the boundary
condition in the usual sense (note that $u\in C^{1}(\overline{\Omega})$). We
say that $u$ is \textit{nontrivial} if $u\not \equiv 0$ and \textit{positive}
if $u>0$ in $\Omega$.

The main feature of this problem is the lack of strong maximum principle
structure, due to the fact that $0<q<1$ and $a$ changes sign. Consequently a
nontrivial solution of $(P_{\alpha})$ is \textit{not} necessarily positive. As
a matter of fact, one may easily find examples where $(P_{\alpha})$ has a
nontrivial solution which is not positive (also known as \textit{dead core}
solution),
%\marginpar{added}
see for instance Remark \ref{2.6} below.
Let us point out that when $a\geq0$ (the definite case) or $q\geq1$ (the
linear and superlinear cases) the strong maximum principle and Hopf's lemma
apply, so in these cases \textit{any} nontrivial solution of $(P_{\alpha})$
belongs to
\[
\mathcal{P}^{\circ}:=\left\{
\begin{array}
[c]{ll}%
\left\{  u\in C_{0}^{1}(\overline{\Omega}):u>0\ \mbox{in $\Omega$},\ \partial
_{\nu}u<0\ \mbox{on $\partial \Omega$}\right\}  &
\mbox{if $\alpha = -\infty$},\medskip\\
\left\{  u\in C^{1}(\overline{\Omega}%
):u>0\ \mbox{on $\overline{\Omega}$}\right\}  &
\mbox{if $\alpha \neq -\infty$}.
\end{array}
\right.
\]

The investigation of $(P_{\alpha})$ in the sublinear case has been carried out
mostly for $\alpha=-\infty$ \cite{BPT1,Br,DS,GK1,GK2,HMV,KRQU16,KRQU3,PT} and
$\alpha=0$ \cite{alama,BPT2,DS,KRQU16,KRQUnodea}. To recall these results, we
consider the conditions
\[
\int_{\Omega}a<0,\leqno{({\bf A.0})}
\]%
\[
%\left\{%
\begin{array}
[c]{l}%
\mbox{$\Omega^{a}_+$ has a finite number of connected components,}
\end{array}
%\right.
\leqno{({\bf A.1})}
\]
where $\Omega_{+}^{a}$ is the open set given by
\[
\Omega_{+}^{a}:=\{x\in\Omega:a(x)>0\}.
\]
We also introduce the positivity set
\[
\mathcal{A}_{\alpha}(a):=\{q\in
(0,1):\mbox{any nontrivial solution of $(P_{\alpha})$ lies in $\mathcal{P}^\circ$}\}.
\]

%\marginpar{: added below}%
\noindent To simplify the notation we write $\mathcal{A}_{\alpha}$ instead of
$\mathcal{A}_{\alpha}(a)$. Note
%\marginpar{{\tiny "No indent" put}}
that $\mathcal{A}_{\alpha}=(0,1)$ whenever $(P_{\alpha})$ has no nontrivial
solution. When $\alpha=0$ (respect. $\alpha=-\infty$) we denote $\mathcal{A}%
_{\alpha}$ by $\mathcal{A}_{N}$ (respect. $\mathcal{A}_{D}$).

We gather now the main results known for $(P_{\alpha})$ in the sublinear case,
which are established in \cite{BPT2}, \cite[Theorem 2.1]{DS}, \cite[Theorems
1.6 and 1.7, Corollary 1.8]{KRQU16}, \cite[Remark 1.1(i)]{KRQUnodea}, and
\cite[Theorem 1.3]{MRS06}:

%Theorem

\begin{theorem}
\label{tin} Let $a\in C(\overline{\Omega})$ be a sign-changing function and
$0<q<1$. Then:

\begin{enumerate}
\item $(P_{-\infty})$ has at least one nontrivial solution.

\item $(P_{0})$ has at least one nontrivial solution under $(A.0)$. Moreover,
if $(P_{0})$ has a positive solution then $(A.0)$ holds.

\item $(P_{\alpha})$ has at most one solution in $\mathcal{P}^{\circ}$ for
$\alpha\in\lbrack-\infty,0]$.

\item Under $(A.1)$ there exists $q_{D}=q_{D}(a)\in\left[  0,1\right)  $ such
that $\mathcal{A}_{D}=(q_{D},1)$. Moreover, if $q\in\mathcal{A}_{D}$ then
$(P_{-\infty})$ has a \textrm{unique} nontrivial solution $u_{D}$, and
$u_{D}\in\mathcal{P}^{\circ}$.

\item Under $(A.0)$ and $(A.1)$ there exists $q_{N}=q_{N}(a)\in\left[
0,1\right)  $ such that $\mathcal{A}_{N}=(q_{N},1)$. Moreover, if
$q\in\mathcal{A}_{N}$ then $(P_{0})$ has a \textrm{unique} nontrivial solution
$u_{N}$, and $u_{N}\in\mathcal{P}^{\circ}$.
\end{enumerate}
\end{theorem}

It is worth pointing out that the uniqueness result in Theorem \ref{tin}(iii)
for the Dirichlet and Neumann problems contrasts with some high multiplicity
results for positive solutions in the superlinear case \cite{BGH,T}. In
Theorem \ref{thm:main}(ii) below we shall prove that for $q\in\mathcal{A}_{N}$
and $\alpha>0$ small $(P_{\alpha})$ has \textit{exactly} two positive
solutions, which shows that a high multiplicity result does not occur in this
situation either.

\begin{figure}[tbh]
\centerline{
\includegraphics[scale=0.14]{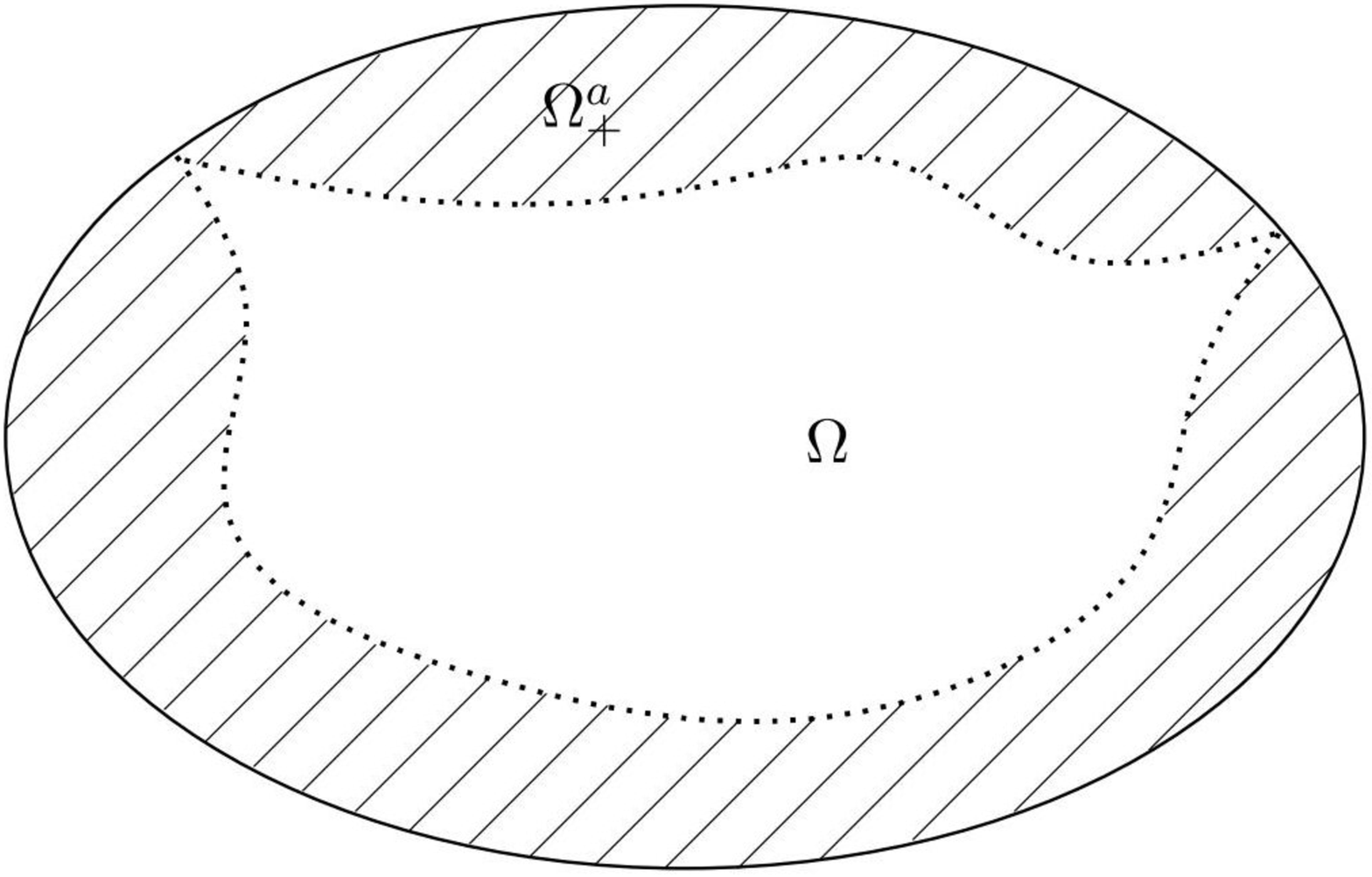}
\hskip0.35cm
\includegraphics[scale=0.14] {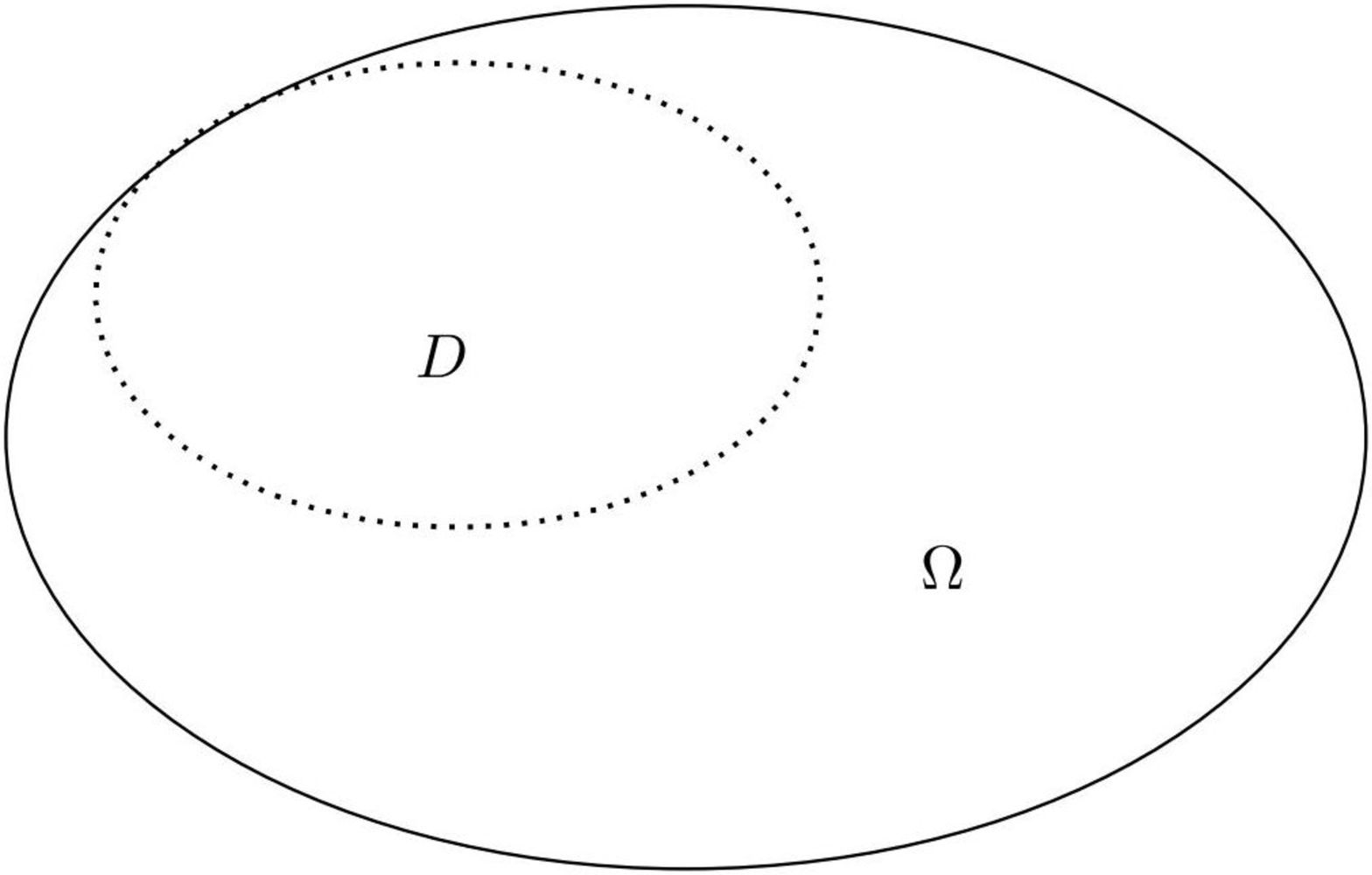}
} \centerline{(i) \hskip6.0cm (ii) }\caption{(i) An example where $(A.2)$
holds; (ii) An example where $(A.3)$ holds.}%
\label{figA2H}%
\end{figure}

In the sequel we state our main results. Some of them shall be established
when $a$ is positive near $\partial\Omega$; more precisely, under the
following assumptions (see Figure \ref{figA2H}):
\[
\partial\Omega\subseteq\partial\Omega_{+}^{a};\leqno{({\bf A.2})}
\]%
\[
0\not \equiv a\geq0\text{ in some smooth domain }D\subset\Omega\text{ such
that }\left\vert \partial D\cap\partial\Omega\right\vert
>0.\leqno{({\bf A.3})}
\]
As in \cite{GMRSL09}, we denote by $\Gamma_{\partial\Omega}$ and
$\Gamma_{\Omega}$ the interior of $\partial D\cap\partial\Omega$ and $\partial
D\cap\Omega$ respectively, and assume that
%$\Gamma_{\partial \Omega}$ and $\Gamma_{\Omega}$ are open;
$\overline{\Gamma_{\partial\Omega}}$, $\overline{\Gamma_{\Omega}}$ are
manifolds with a common $N-2$ dimensional boundary $\Gamma^{\prime}%
:=\overline{\Gamma_{\partial\Omega}}\cap\overline{\Gamma_{\Omega}}$, and
$\partial D=\Gamma_{\partial\Omega}\cup\Gamma^{\prime}\cup\Gamma_{\Omega}$.

The main role of $(A.2)$ is to ensure that any nontrivial solution of
$(P_{\alpha})$ satisfies $u\not \equiv 0$ in $\Omega_{+}^{a}$ for any
$\alpha>0$, cf. Lemma \ref{lpos}. As for $(A.3)$, it shall provide us with
\textit{a priori} bounds on $\alpha>0$ for the existence of solutions in
$\mathcal{P}^{\circ}$, cf. Propositions \ref{alpha2} and \ref{prop:nbddalph}.
Let us note that $(A.2)$ holds if $a>0$ at every point on $\partial\Omega$;
nevertheless, $(A.2)$ may still be true if $a$ vanishes (somewhere or
everywhere) on $\partial\Omega$.

%\marginpar{\scriptsize need modify}
%\begin{figure}[tbh]
%\begin{center}
%\includegraphics[scale=0.15]{fig18_1202H.eps}
%\end{center}
%\caption{An example where $(H)$ holds.}%
%\label{figsubdD02}%
%\end{figure}

%\begin{figure}[H]
%\begin{figure}[tbh]
%\begin{center}
%\includegraphics[scale=0.18]{fig18_1203A2.eps}
%\end{center}
%\caption{An example where $(A.2)$ holds.}
%\label{figA23}
%\end{figure}

%\marginpar{\scriptsize Fig.\ref{figA23} modified}

%\begin{figure}[tbh]
%\centerline{
%\includegraphics[scale=0.14]{fig18_1203A2.eps}
%\hskip0.35cm
%\includegraphics[scale=0.14] {fig18_1202A3.eps}
%} \centerline{(i) \hskip6.0cm (ii) }\caption{(i) An example where $(A.2)$
%holds; (ii) An example where $(A.3)$ holds.}%
%\label{figA23}%
%\end{figure}

We start by showing that $(P_{\alpha})$ inherits the positivity property from
the Dirichlet problem (i.e. for $q\in\mathcal{A}_{D}$) up to a certain
$\alpha_{p}(a)>0$:

\begin{theorem}
[Positivity]\label{tpos} Assume $(A.1)$. Then there exists $\alpha_{p}%
=\alpha_{p}(a)>0$ such that any nontrivial solution of $(P_{\alpha})$ belongs
to $\mathcal{P}^{\circ}$ for every $\alpha<\alpha_{p}$ and $q \in
\mathcal{A}_{D}$. Moreover, $\alpha_{p}=\infty$ if $(A.2)$ holds.
\end{theorem}

In view of the above theorem, we shall deal with $q\in\mathcal{A}_{D}$ in most
of our results. We proceed with the description of the solution set of
$(P_{\alpha})$ for $\alpha<0$. This case turns out to be similar to the
Dirichlet one, as long as existence and uniqueness of a nontrivial solution
are concerned. As a matter of fact, when $\alpha<0$ we shall see that $(A.0)$
is \textit{not} necessary for the existence of a positive solution, unlike in
the case $\alpha\geq0$ (for the Neumann problem see \cite[Lemma 2.1]{BPT2},
which can be easily extended to $\alpha>0$).

%Theorem

\begin{figure}[tbh]
\centerline{
\includegraphics[scale=0.18]{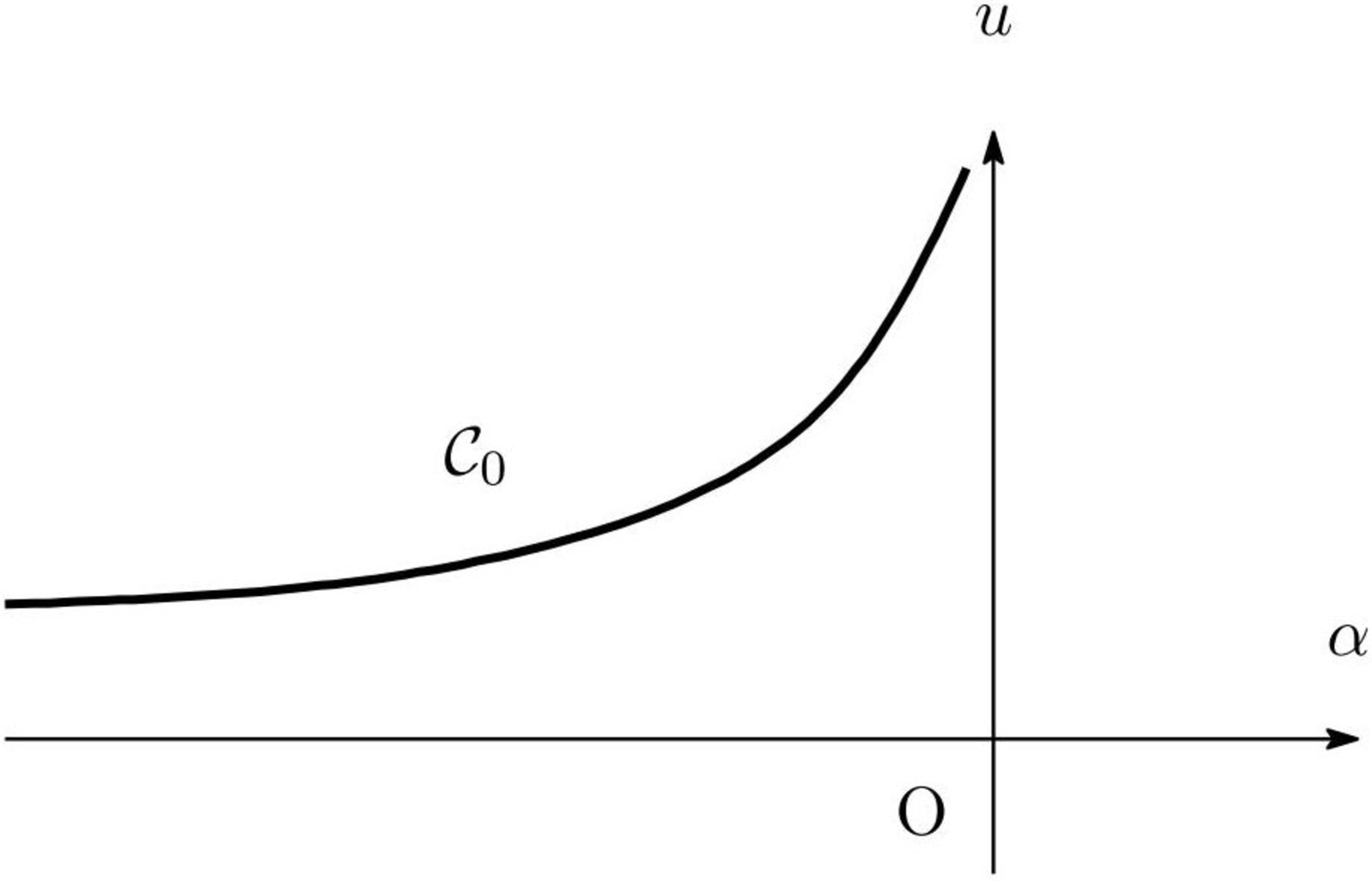} \hskip0.35cm
\includegraphics[scale=0.18] {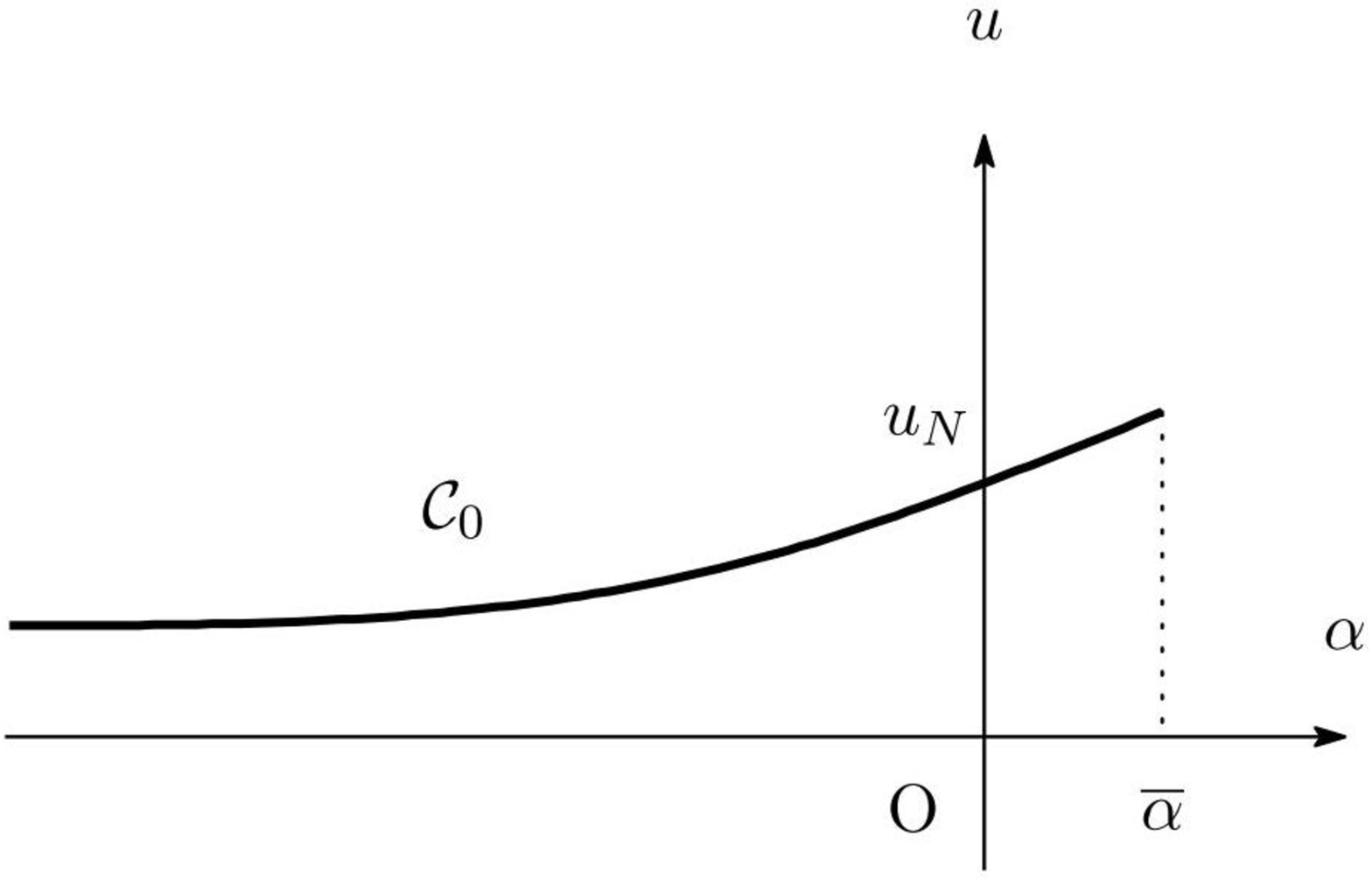}
} \centerline{(i) \hskip5.0cm (ii) }\caption{The positive solutions curve
$\mathcal{C}_{0}$ for $q \in\mathcal{A}_{D}$ in the cases (i) $\int_{\Omega
}a\geq0$ and (ii) $\int_{\Omega}a <0$.}%
\label{fig:do1}%
\end{figure}

\begin{theorem}
[A curve of positive solutions for $\alpha<0$]\label{thm:curve} Assume $(A.1)$
and $q\in\mathcal{A}_{D}$. Then $(P_{\alpha})$ has a unique nontrivial
solution $u_{\alpha}$ for each $\alpha<0$, and $u_{\alpha}\in\mathcal{P}%
^{\circ}$. Moreover, the mapping $\mathcal{C}_{0}:\alpha\mapsto u_{\alpha}$ is
$C^{\infty}$ from $(-\infty,0)$ into $W^{2,r}(\Omega)$, increasing (i.e.
$u_{\alpha}<u_{\beta}$ on $\overline{\Omega}$ for $\alpha<\beta<0$), and
$u_{\alpha}\rightarrow u_{D}$ in $H^{1}(\Omega)$ as $\alpha\rightarrow-\infty
$. Finally, as $\alpha\rightarrow0^{-}$ we have the following alternative:

\begin{enumerate}
\item Assume that $\left(  A.0\right)  $ does
%\marginpar{$H_{1}\left(\Omega\right)  $ {\tiny added to the norm, since here we do not yet defined}$\left\Vert {}\right\Vert $}
not hold. Then $\displaystyle\min_{x\in\overline{\Omega}}u_{\alpha
}(x)\rightarrow\infty$ as $\alpha\rightarrow0^{-}$ (see Figure \ref{fig:do1}%
(i)). In particular, $u_{\alpha}$ approaches a spatially homogeneous
distribution on $\overline{\Omega}$. Moreover, $(P_{\alpha})$ has no solution
$u$ such that $u\not \equiv 0$ in $\Omega_{+}^{a}$ for $\alpha\geq0$.

\item Assume that $\left(  A.0\right)  $ holds. Then $\mathcal{C}_{0}$ can be
extended to $(-\infty,\overline{\alpha})$, for some $\overline{\alpha}>0$, so
that $u_{0}=u_{N}$ and $u_{\alpha}\in\mathcal{P}^{\circ}$ solves $(P_{\alpha
})$ for $\alpha\in(0,\overline{\alpha})$. Moreover, the mapping $\alpha\mapsto
u_{\alpha}$ is increasing in $(-\infty,\overline{\alpha})$ and unique in the
following sense: if $u_{n}$ solves $(P_{\alpha_{n}})$ with $\alpha
_{n}\rightarrow0^{+}$ and $u_{n}\rightarrow u_{N}$ in $C^{1}(\overline{\Omega
})$, then, for $n$ large enough, $u_{n}=u_{\alpha}$ for some $\alpha
\in(0,\overline{\alpha})$ (see Figure \ref{fig:do1}(ii)).
\end{enumerate}
\end{theorem}

\begin{rem}
\strut

\begin{enumerate}
\item Let $0<q<1$. Under $(A.1)$ there exists $\delta=\delta(q,a^{+})>0$ such
that any nontrivial solution of $(P_{-\infty})$ lies in $\mathcal{P}^{\circ}$
if $0<\Vert a^{-}\Vert_{C(\overline{\Omega})}<\delta$, cf. \cite[Theorem
1.1]{KRQU16}. One may easily see from its proof that Theorem \ref{thm:curve}
also holds if we assume $0<\Vert a^{-}\Vert_{C(\overline{\Omega})}<\delta$,
instead of $q\in\mathcal{A}_{D}$.

\item A 'bifurcation from infinity' scenario as in Theorem \ref{thm:curve}(i)
also occurs under $(A.0)$, now with $\alpha\rightarrow0^{+}$ (see Theorem
\ref{thm:main}(ii-c)).
\end{enumerate}
\end{rem}

%	%	% \begin{figure}[H]
%\begin{figure}[tbh]
%\begin{center}
%\includegraphics[scale=0.18]{fig18_1004curve2.eps}
%\end{center}
%\caption{The positive solutions curve $\mathcal{C}_{0}$ when $\int_{\Omega}a
%\geq0$ and $q \in\mathcal{A}_{D}$.}%
%\label{fig18_1004curve2}%
%\end{figure}

%\marginpar{The caption modified in Fig. 1 and 2}
%	%	% \begin{figure}[H]
%\begin{figure}[tbh]
%\begin{center}
%\includegraphics[scale=0.18]{fig18_0915curve.eps}
%\end{center}
%\caption{The positive solutions curve $\mathcal{C}_{0}$ when $\int_{\Omega}a <
%0$ and $q \in\mathcal{A}_{D}$.}%
%\label{fig18_0915curve}%
%\end{figure}

Differently from the case $\alpha\leq0$, we shall see that when $\alpha>0$ is
small enough $(P_{\alpha})$ may admit \textit{multiple} solutions in
$\mathcal{P}^{\circ}$. To this end, we set
\begin{equation}
w:=\alpha^{\frac{1}{1-q}}u \label{changeofv}%
\end{equation}
and transform $(P_{\alpha})$ into
\[%
\begin{cases}
-\Delta w=\alpha a(x)w^{q} & \mbox{in $\Omega$},\\
w\geq0 & \mbox{in $\Omega$},\\
\partial_{\nu}w=\alpha w & \mbox{on $\partial \Omega$}.
\end{cases}
\leqno{(R_{\alpha})}
\]
%Now, setting $w=\alpha^{\frac{1}{1-q}}v$ we may transform $(P_{\alpha})$ into
%\[%
%\begin{cases}
%-\Delta w=\alpha a(x)w^{q} & \mbox{in $\Omega$},\\
%\partial_{\nu} w=\alpha w & \mbox{on $\partial \Omega$}.
%\end{cases}
%\leqno{(R_{\alpha})}
%\]
%\marginpar{removed the sentence about Chabrowski-Tintarev}
We shall treat this problem via a bifurcation approach, looking at $\alpha$ as
a bifurcation parameter. It turns out that $(R_{\alpha})$ is easier to handle
(in comparison with $(P_{\alpha})$), providing us with a more accurate
description of the solutions set of $(P_{\alpha})$ for $\alpha>0$ small.
Indeed, note that $(R_{\alpha})$ has two solutions lines, namely:
\begin{equation}%
\begin{cases}
\Gamma_{0}:=\{(\alpha,w):\alpha>0,\ w=0\}, & \\
\Gamma_{1}:=\{(\alpha,w):\alpha=0,\ \mbox{$w$ is
a nonnegative constant}\}. &
\end{cases}
\label{G0:00}%
\end{equation}
Under $(A.0)$, let us put
\begin{equation}
c_{a}:=\left(  \frac{-\int_{\Omega}a}{|\partial\Omega|}\right)  ^{\frac
{1}{1-q}}. \label{def:ca}%
\end{equation}
In \cite[Section 7]{CT14} Chabrowski and Tintarev proved, by variational
methods, that under $(A.0)$ this problem has at least two nontrivial solutions
$w_{1,\alpha},w_{2,\alpha}$ such that $w_{1,\alpha}<w_{2,\alpha}$ on
$\overline{\Omega}$ for $\alpha>0$ small enough. Moreover, they also provided
the following asymptotic profiles of $w_{1,\alpha},w_{2,\alpha}$ as
$\alpha\rightarrow0^{+}$:
\begin{equation}
w_{2,\alpha}\rightarrow c_{a}\quad\mbox{and}\quad w_{1,\alpha}\rightarrow
0\ \mbox{in}\ H^{1}(\Omega)\quad\mbox{as}\quad\alpha\rightarrow0^{+},
\label{w1w2}%
\end{equation}
%where $c_{a}$ is defined by \eqref{def:ca},
and every sequence $\alpha_{n}\rightarrow0$ has a subsequence (still denoted
by the same notation) satisfying
\begin{equation}
\alpha_{n}^{-\frac{1}{1-q}}w_{1,\alpha_{n}}\rightarrow u_{0}%
\mbox{ in $H^{1}(\Omega)$}, \label{w1ap}%
\end{equation}
where $u_{0}$ is a nontrivial solution of $(P_{0})$.

We shall complement \eqref{w1w2} and \eqref{w1ap} in two ways, proving the following:

\begin{itemize}
\item[(I)] an \textit{exact} multiplicity result for $q\in\mathcal{A}_{N}$,
namely: $w_{1,\alpha},w_{2,\alpha}$ are the only nontrivial solutions of
$(R_{\alpha})$ if $\alpha>0$ is small enough, and $w_{1,\alpha},w_{2,\alpha
}\in\mathcal{P}^{\circ}$ (Theorem \ref{prop:exactm});

\item[(II)] the existence of a subcontinuum of solutions of $(R_{\alpha})$ for
$\alpha>0$ small, connecting $(0,0)$ to $(0,c_{a})$ (Theorem \ref{maint}, see
also Remark \ref{4.5}).
\end{itemize}

\begin{figure}[tbh]
\centerline{
\includegraphics[scale=0.175]{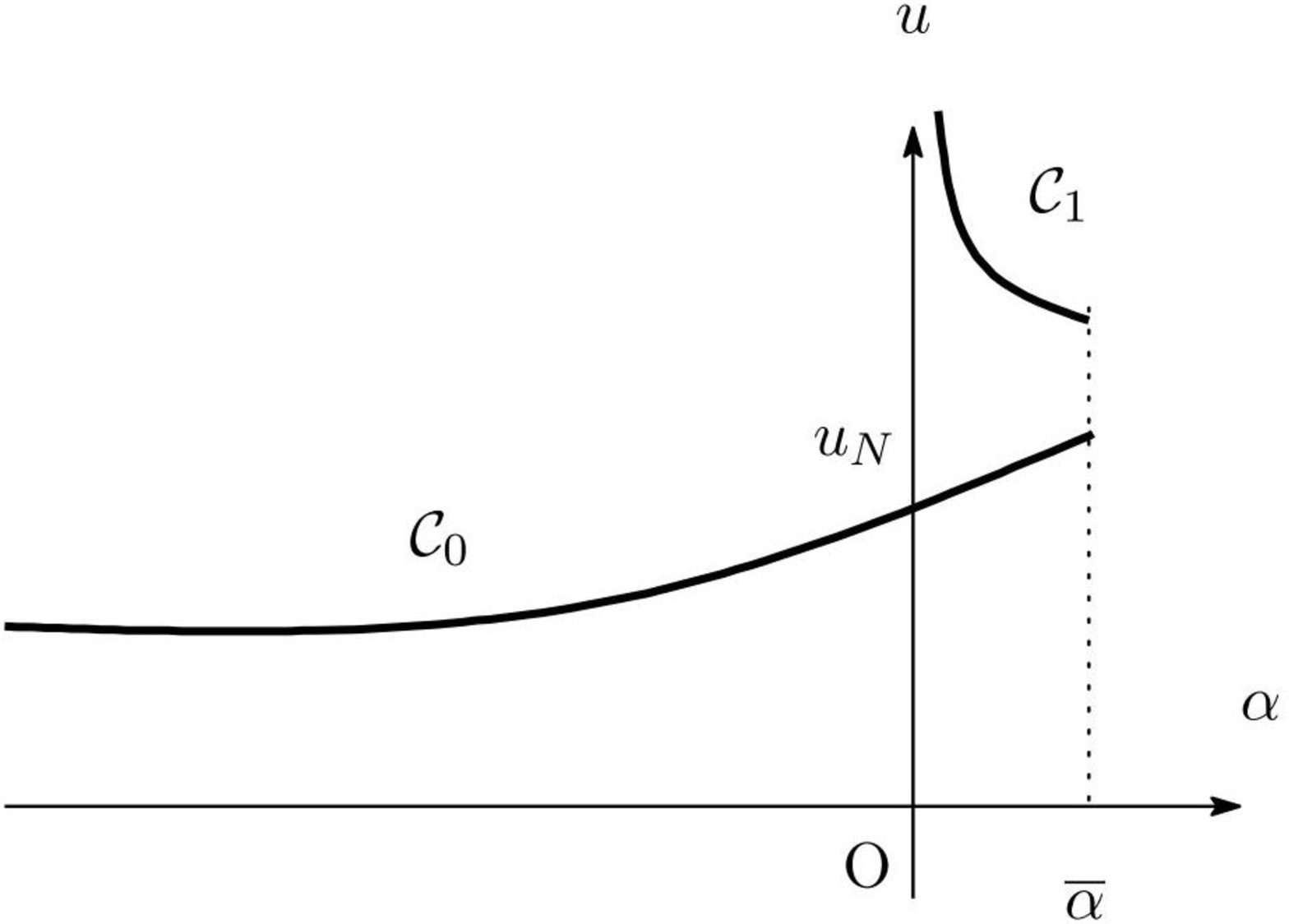} \hskip0.35cm
\includegraphics[scale=0.175] {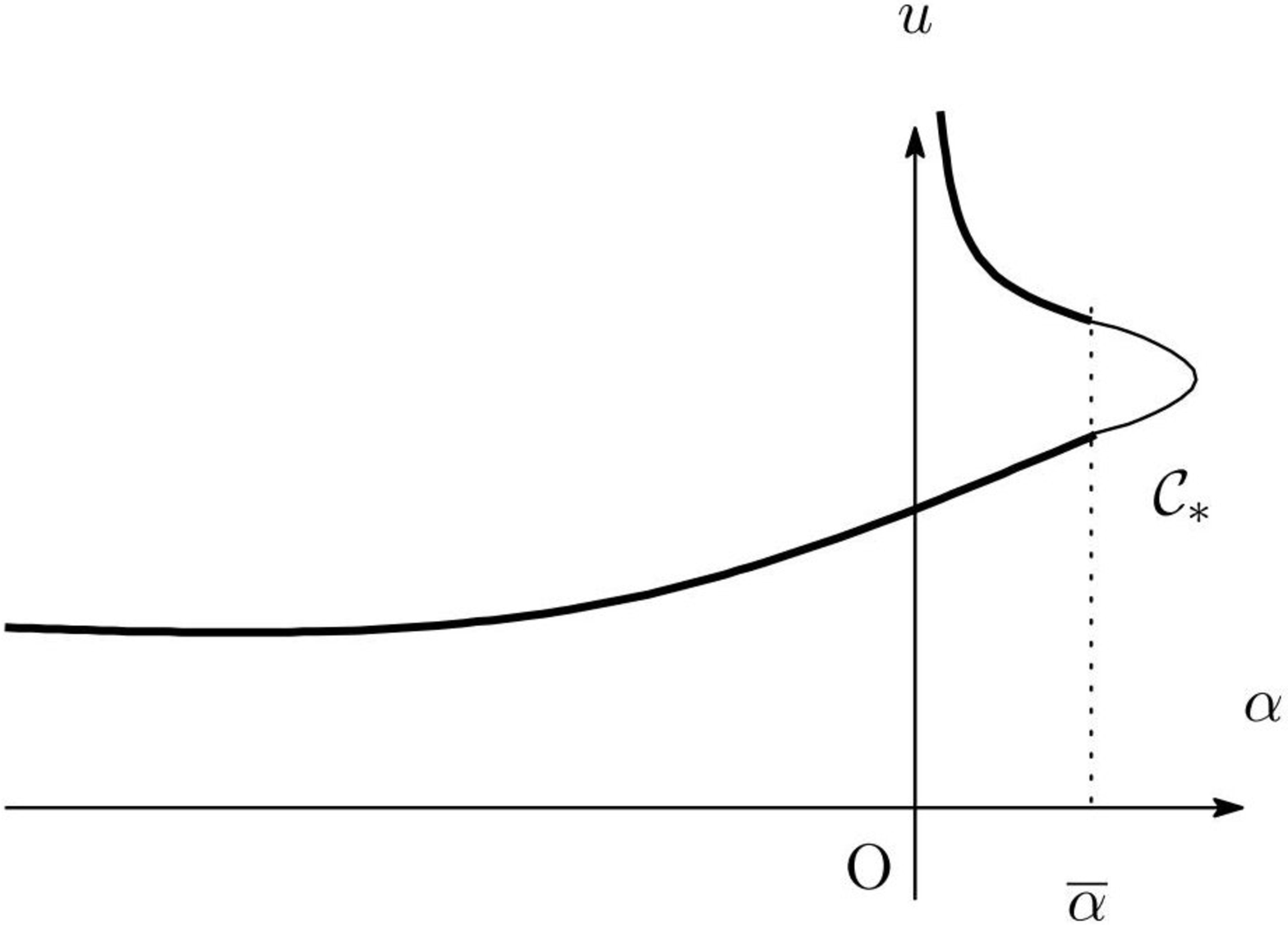}
} \centerline{(i) \hskip5.0cm (ii) }\caption{(i) The $C^{\infty}$ curves
$\mathcal{C}_{0}$ and $\mathcal{C}_{1}$, (ii) The component $\mathcal{C}%
_{\ast}$.}%
\label{fig:do2}%
\end{figure}

These results, combined with Theorems \ref{tpos} and \ref{thm:curve}, provide
a global description (with respect to $\alpha$) of the solutions set of
$(P_{\alpha})$ for $q\in\mathcal{A}_{D}$:

\begin{theorem}
\label{thm:main}
%\marginpar{modified, statement and proof}
Assume $(A.0),(A.1)$, and $q\in\mathcal{A}_{D}$. Then the following assertions hold:

\begin{enumerate}
\item (Existence and nonexistence) Let
\begin{equation}
\alpha_{s}:=\sup\{\alpha\in\mathbb{R}%
:\mbox{$(P_\alpha)$ has a solution in $\mathcal{P}^\circ$}\}. \label{as}%
\end{equation}
Then $\alpha_{s}\in\left(  0,\infty\right)  $, i.e. $(P_{\alpha})$ has at
least one solution in $\mathcal{P}^{\circ}$ for $\alpha>0$ small and no such
solution for $\alpha>0$ large. In addition, if $(A.2)$ holds then $(P_{\alpha
})$ has at least one solution in $\mathcal{P}^{\circ}$ for every $\alpha
\leq\alpha_{s}$.

\item (Exact multiplicity and limiting behavior) There exists $\overline
{\alpha}\in(0,\alpha_{s}]$ such that $(P_{\alpha})$ has exactly one nontrivial
solution $u_{1,\alpha}$ for $\alpha\in(-\infty,0]$, and exactly two nontrivial
solutions $u_{1,\alpha},u_{2,\alpha}$ for $\alpha\in(0,\overline{\alpha})$.
Moreover $u_{1,\alpha}, u_{2,\alpha}\in\mathcal{P}^{\circ}$ and these ones
satisfy (see Figure \ref{fig:do2}(i)):

\begin{enumerate}
\item $\mathcal{C}_{0}:\alpha\mapsto u_{1,\alpha}$ is a $C^{\infty}$
increasing mapping from $(-\infty,\overline{\alpha})$ into $W^{2,r}(\Omega)$
(as in Theorem \ref{thm:curve}(ii)), and $\mathcal{C}_{1}:\alpha\mapsto
u_{2,\alpha}$ is a $C^{\infty}$ mapping from $(0,\overline{\alpha})$ into
$W^{2,r}(\Omega)$.

\item $u_{1,\alpha}\rightarrow u_{D}$ in $H^{1}(\Omega)$ as $\alpha
\rightarrow-\infty$.

\item $u_{2,\alpha}\sim\alpha^{-\frac{1}{1-q}}c_{a}$ as $\alpha\rightarrow
0^{+}$, i.e. $\alpha^{\frac{1}{1-q}}u_{2,\alpha}\rightarrow c_{a}$ in
$C^{1}(\overline{\Omega})$ (implying $\displaystyle \min_{\overline{\Omega}%
}u_{2,\alpha}\rightarrow\infty$) as $\alpha\rightarrow0^{+}$, where $c_{a}$ is
given by \eqref{def:ca}.
\end{enumerate}

\item (Existence of a component) Assume in addition $(A.2)$ and $\left(
A.3\right)  $. Then $(P_{\alpha})$ possesses a component $\mathcal{C}_{\ast}$
(\textrm{i.e. a maximal closed, connected subset} in $\mathbb{R}\times
C^{1}(\overline{\Omega})$) of solutions in $\mathcal{P}^{\circ}$ that contains
$\mathcal{C}_{0}$ and $\mathcal{C}_{1}$. In addition,
\[
\mathcal{C}_{\ast}\cap\left\{  (\alpha,0):\alpha\in\mathbb{R}\right\}
=\emptyset,
\]
and
\[
\mathcal{C}_{\ast}\cap\{(\alpha,\infty):\alpha\in\mathbb{R}\}=\{(0,\infty)\},
\]
i.e. $\mathcal{C}_{\ast}$ does not meet the trivial solution at any $\alpha
\in\mathbb{R}$ and $\mathcal{C}_{\ast}$ bifurcates from infinity only at
$\alpha=0$, see Figure \ref{fig:do2}(ii).
%\marginpar{modified}

%\begin{itemize}
%\item $\sup\left\{  \alpha\in\mathbb{R}:(\alpha,u)\in\mathcal{C}_{\ast}\right\}  <\infty$.
%\item
%$\mathcal{C}_{\ast}\cap\{(\alpha,0),(\alpha,\infty):\alpha\geq\overline{\alpha}\}=\emptyset$.
%\end{itemize}

\end{enumerate}
\end{theorem}

%	%	% \begin{figure}[H]
%\begin{figure}[tbh]
%\begin{center}
%\includegraphics[scale=0.2]{fig18_0918C1C2.eps}
%\end{center}
%\caption{The $C^{\infty}$ curves $\mathcal{C}_{1}$ and $\mathcal{C}_{2}$.}%
%\label{fig18_0918C1C2}%
%\end{figure}
%	%	% \begin{figure}[H]
%\begin{figure}[tbh]
%\begin{center}
%\includegraphics[scale=0.2]{fig18_0918C1C2C0.eps}
%\end{center}
%\caption{The component $\mathcal{C}_{0}$.}%
%\label{fig18_0918C0}%
%\end{figure}

\begin{rem}
\strut\label{N}

\begin{enumerate}
\item One may easily check that all the assertions for $\alpha\geq0$ in
Theorem \ref{thm:main} still hold if we take $q\in\mathcal{A}_{N}$ in place of
$q\in\mathcal{A}_{D}$ (let us note that $\mathcal{A}_{N}\supseteq
\mathcal{A}_{D}$ by Proposition \ref{cor:Aal:ordered} below).

\item Some lower and upper bounds on $\alpha_{s}$ are given in Corollary
\ref{cal}. Moreover, we
%\marginpar{modified}
shall provide (finite) upper bounds for $\alpha_{s}$ for \textit{every}
$q\in\lbrack0,1)$, see Proposition \ref{alpha2} below.

\item The approach to obtain the solution $u_{2,\alpha}$ from Theorem
\ref{thm:main}(ii) applies to any $q\in\lbrack0,1)$. Thus, for $q$ close to
$0$ (including $q=0$), $(P_{\alpha})$ still has, under $(A.0)$, a solution in
$\mathcal{P}^{\circ}$ for $\alpha>0$ small, see Remark \ref{rpq}. We note that
when $q=0$ and $\int_{\Omega}a<0$ there are no solutions of $\left(
P_{\alpha}\right)  $ in $\mathcal{P}^{\circ}$ for $\alpha\in\left[
-\infty,0\right]  $, since in this case solutions $u$ satisfy $-\int_{\Omega
}a=\alpha\int_{\partial\Omega}u$ if $\alpha\neq-\infty$ and $-\int_{\Omega
}a=\int_{\partial\Omega}\partial_{\nu}u$ if $\alpha=-\infty$.
%\marginpar{added}

\item Theorems \ref{tpos}, \ref{thm:curve} and \ref{thm:main} hold more
generally for $a\in L^{r}(\Omega)$, with $r>N$. In this case, we set
$\Omega_{+}^{a}$ as the largest open subset of $\Omega$ in which $a>0$ a.e.,
and we add to $(A.1)$ the condition $|(\text{supp}\,a^{+})\setminus\Omega
_{+}^{a}|=0$, where supp is the support in the measurable sense.
\end{enumerate}
\end{rem}

To the best of our knowledge, \textit{exact multiplicity} results are not
commonly seen in the literature, specially for indefinite type problems such
as $(P_{\alpha})$. We refer to \cite[Section 3]{KO} for a result of this kind
with $N=1$ and a superlinear nonlinearity.
%\marginpar{added}
Let us add that some multiplicity results for $(P_{-\infty})$ and $(P_{0})$
are given in \cite[Section 2]{BPT1} and \cite[Section 4]{BPT2}, \cite[Theorem
1.1]{alama} respectively.

Finally, although we are mainly focused on $q\in(0,1)$, we shall see that when
$q=0$ and $\alpha>0$ many interesting questions arise. Some of them are
treated in this article, whereas some other ones are left to a forthcoming paper.

The rest of the article is organized as follows. In Section \ref{secneg} we
mainly analyze the case $\alpha\leq0$ and prove Theorems \ref{tpos} and
\ref{thm:curve}. Section \ref{secpos} is mostly devoted to $(R_{\alpha})$ with
$\alpha>0$, where we investigate qualitative properties of the solutions set
and prove an exact multiplicity result employing the change of variables
\eqref{changeofv}. Lastly, Section \ref{sec:topoba} provides a topological
bifurcation approach of $(R_{\alpha})$ and the proof of Theorem \ref{thm:main}.

\subsubsection*{Notation}

\begin{itemize}
\setlength{\itemsep}{0.0cm}

\item For any $f\in L^{1}(\Omega)$ the integral $\int_{\Omega}f$ is considered
with respect to the Lebesgue measure, whereas for any $g\in L^{1}%
(\partial\Omega)$ the integral $\int_{\partial\Omega}g$ is considered with
respect to the surface measure.

\item The usual norm of $H^{1}(\Omega)$ is denoted by $\Vert\cdot\Vert$, i.e.
$\Vert u\Vert=\left(  \int_{\Omega}|\nabla u|^{2}+u^{2}\right)  ^{\frac{1}{2}%
}$. For $r\geq1$ the Lebesgue norm in $L^{r}(\Omega)$ will be denoted by
$\Vert\cdot\Vert_{r}$.
%\marginpar{added}

\item The weak convergence is denoted by $\rightharpoonup$.

\item The positive and negative parts of a function $u$ are defined by
$u^{\pm}:=\max\{\pm u,0\}$.

\item $|\cdot|$ stands for both the Lebesgue measure and the surface measure.

\item The characteristic function of a set $A\subset\mathbb{R}^{N}$ is denoted
by $\chi_{A}$.
%\marginpar{added}
%\marginpar{modified}
%\marginpar{added}
%\item If $U \subset \R^N$ then we denote its interior by int $U$ and its closure by $\overline{U}$. The characteristic function of $U$ is denoted by $\chi_U$.
%\item We denote by $B(0,R)$ the ball of radius $R$ centered at $0$ in $H_0^1(\Omega)$.

\end{itemize}

\section{Proof of Theorems \ref{tpos} and \ref{thm:curve}}

\label{secneg}

We split the proofs of Theorems \ref{tpos} and \ref{thm:curve} into several
results. The first one is a direct consequence of Lemma \ref{lpos} and
Proposition \ref{cor:Aal:ordered}, whereas the second one follows from
Propositions \ref{alinf}, \ref{pc} and \ref{last}.

We start by proving that nontrivial solutions of $(P_{\alpha})$ are positive
in some component of $\Omega^{a}_{+}$ as long as $\alpha$ is less than
\[
\alpha_{p}=\alpha_{p}(a):=\inf\left\{  \int_{\Omega}|\nabla v|^{2}: v \in
H^{1}(\Omega), v \equiv0 \text{ in } \Omega^{a}_{+}, \int_{\partial\Omega}
v^{2}=1\right\}  .
\]
Note that $\alpha_{p}$ depends on $a$ but not on $q$.

\begin{lem}
\label{lpos} \strut

\begin{enumerate}
\item We have $\alpha_{p}>0$. Moreover $\alpha_{p}=\infty$ if, and only if,
$(A.2)$ holds.

\item We have $u\not \equiv 0$ in $\Omega_{+}^{a}$ for any nontrivial solution
$u$ of $(P_{\alpha})$ and for any $\alpha<\alpha_{p}$ and $q\in\lbrack0,1)$.
\end{enumerate}
\end{lem}

\noindent\textit{Proof.}

\begin{enumerate}
\item First of all, one may easily show that this infimum is achieved whenever
it is finite, and consequently that it is positive, since no constant function
satisfies the constraints simultaneously. Now, if $(A.2)$ holds then there is
no $v$ satisfying $v \equiv0$ in $\Omega_{a}^{+}$ and $\int_{\partial\Omega}
v^{2}=1$, so that $\alpha_{p}=\infty$. Finally, if $(A.2)$ does not hold then
we may find some ball $B$ around some $x_{0} \in\partial\Omega$ such that $a
\leq0$ in $B \cap\Omega$. We may then build some $v \in H^{1}(\Omega)$
supported in $B \cap\Omega$ and such that $\int_{\partial\Omega} v^{2}=1$.
Thus $v$ is admissible for $\alpha_{p}$, and consequently $\alpha_{p}<\infty$.

\item Let $\alpha<\alpha_{p}$ and $u$ be a nontrivial solution of $(P_{\alpha
})$. If $u\equiv0$ in $\Omega_{a}^{+}$ then we have
\[
0<\int_{\Omega}|\nabla u|^{2}=\int_{\Omega}a(x)u^{q+1}+\alpha\int
_{\partial\Omega}u^{2}\leq\alpha\int_{\partial\Omega}u^{2},
\]
so that $\alpha\geq(\int_{\Omega}|\nabla u|^{2})(\int_{\partial\Omega}%
u^{2})^{-1}$. Consequently $\alpha\geq\alpha_{p}$, which contradicts our
assumption. \qed\newline
\end{enumerate}

\begin{rem}
Assume that $\Omega_{0,+}^{a}:=\{x\in\Omega:a(x)\geq0\}$ is connected and
smooth. Then $\alpha_{p}$ can be reset as
\[
\alpha_{p}=\alpha_{p}(a):=\inf\left\{  \int_{\Omega}|\nabla v|^{2}: v \in
H^{1}(\Omega), v \equiv0 \text{ on } \Omega_{0,+}^{a}, \int_{\partial\Omega}
v^{2}=1\right\}  .
\]
In this case, Lemma \ref{lpos}(i) holds with $(A.2)$ formulated now as
$\partial\Omega\subseteq\partial\Omega_{0,+}^{a}$. Moreover, one can repeat
the proof of Lemma \ref{lpos}(ii) to show that $u \not \equiv 0$ on
$\Omega^{a}_{0,+}$ for any nontrivial solution $u$ of $(P_{\alpha})$ and any
$\alpha<\alpha_{p}$. Since $\Omega_{0,+}^{a}$ is smooth and connected, the
strong maximum principle yields $u >0$ in $\Omega_{+}^{a}$. Note that this new
value $\alpha_{p}$ is larger than the original one.
\end{rem}

\begin{prop}
[Monotonicity of $\mathcal{A}_{\alpha}$]\label{cor:Aal:ordered}
%\marginpar{modified statement and proof}
We have $\mathcal{A}_{\alpha}\subseteq\mathcal{A}_{\beta}$ for $-\infty
\leq\alpha<\beta<\alpha_{p}$.
\end{prop}

\noindent\textit{Proof.} First we consider $\alpha=-\infty$. Let
%\marginpar{modified}
$q\in\mathcal{A}_{D}$ and $u$ be a nontrivial solution of $(P_{\beta})$. Since
$u\geq0$ on $\overline{\Omega}$, we see that $u$ is a supersolution of
$(P_{-\infty})$. Moreover, by Lemma \ref{lpos} we know that $u\not \equiv 0$
in $\Omega_{+}^{a}$.
%Indeed, if not then
%\[
%0<\int_{\Omega}|\nabla u|^{2}=\int_{\Omega}a(x)u^{q+1}+\beta\int %_{\partial\Omega}u^{2}\leq0.
%\]
It follows that there exist a ball $B\subset\Omega_{+}^{a}$ and a constant
$c>0$ such that $u>c$ in $B$. It is then possible to provide a subsolution $z$
of $(P_{-\infty})$ such that $z\not \equiv 0$, $\mathrm{supp\,}z\Subset B$,
and $z\leq u$
%\marginpar{\textit{added, and put above "provide" place of "construct" to avoid repetition}}
(see e.g. the construction in \cite[Lemma 2.3(ii)]{BPT1}). By the sub and
supersolution method, we find a nontrivial solution $u_{1}$ of $(P_{-\infty})$
such that $z\leq u_{1}\leq u$ on $\overline{\Omega}$. Since $q\in
\mathcal{A}_{D}$, we have $u_{1}\in\mathcal{P}^{\circ}$, so $u\geq u_{1}>0$ in
$\Omega$ and $\partial_{\nu}u_{1}<0$ on $\partial\Omega$. We claim that $u>0$
on $\partial\Omega$. Indeed, otherwise we have $u=\partial_{\nu}u=0$ somewhere
on $\partial\Omega$. But since $u_{1}=0>\partial_{\nu}u_{1}$ on $\partial
\Omega$, this contradicts the assertion $u_{1}\leq u$ in $\Omega$. Hence $u>0$
on $\overline{\Omega}$, which shows that $q\in\mathcal{A}_{\beta}$.

Let now $\alpha>-\infty$. Take
%\marginpar{modified}
$q\in\mathcal{A}_{\alpha}$ and $u$ a nontrivial solution of $(P_{\beta})$.
Then, arguing as in the previous case, we find by the sub and supersolution
method a nontrivial solution $u_{1}$ of $(P_{\alpha})$ such that $u_{1}\leq u$
on $\overline{\Omega}$. Since $q\in\mathcal{A}_{\alpha}$, it follows that
$u\geq u_{1}>0$ on $\overline{\Omega}$, which shows that $q\in\mathcal{A}%
_{\beta}$. \qed\newline

Next we deal with
\begin{equation}
\tilde{\alpha}=\tilde{\alpha}(q,a):=\inf\left\{  \int_{\Omega}|\nabla
v|^{2}:v\in H^{1}(\Omega),\int_{\Omega}a(x)|v|^{q+1}\geq0,\int_{\partial
\Omega}v^{2}=1\right\}  . \label{def:tilal}%
\end{equation}

One may easily show that this infimum is achieved. Note also that
$\tilde{\alpha}\geq0$ and that $\tilde{\alpha}>0$ if, and only if, $(A.0)$
holds. Lastly, one may show that $\displaystyle \liminf_{q \to1^{-}}
\tilde{\alpha}(q,a) \geq\tilde{\alpha}(1,a)$, so that $\tilde{\alpha}(q,a)$
stays away from zero for $q$ close to $1$ if $(A.0)$ holds.

\begin{prop}
[Existence of a solution in $\mathcal{P}^{\circ}$]\label{alinf}$(P_{\alpha})$
has at least one solution $u_{\alpha}$ such that $u_{\alpha}\not \equiv 0$ in
$\Omega_{a}^{+}$ for every $\alpha<\tilde{\alpha}$. In addition:

\begin{enumerate}
\item Assume $(A.1)$ and $q\in\mathcal{A}_{D}$. Then $u_{\alpha}
\in\mathcal{P}^{\circ}$ and $u_{\alpha}$ is the unique nontrivial solution of
$(P_{\alpha})$ for $\alpha<0$.

\item Assume $(A.0)$, $(A.1)$ and $q\in\mathcal{A}_{N}$. Then $u_{\alpha}%
\in\mathcal{P}^{\circ}$ for $0<\alpha<\tilde{\alpha}$.
\end{enumerate}
\end{prop}

\noindent\textit{Proof.} Let
\begin{equation}
\mu(\alpha):=\inf\left\{  \int_{\Omega}|\nabla u|^{2}-\alpha\int
_{\partial\Omega}u^{2}:u\in H^{1}(\Omega),\int_{\Omega}a(x)|u|^{q+1}%
=1\right\}  \label{mua}%
\end{equation}
We claim that $\mu(\alpha)$ is finite if $\alpha<\tilde{\alpha}$. Indeed,
assume by contradiction that $u_{n}$ satisfies $\int_{\Omega}a(x)|u_{n}%
|^{q+1}=1$ and $\int_{\Omega}|\nabla u_{n}|^{2}-\alpha\int_{\partial\Omega
}u_{n}^{2}\rightarrow-\infty$. In particular, we have $\Vert u_{n}%
\Vert\rightarrow\infty$. We set $v_{n}:=\frac{u_{n}}{\Vert u_{n}\Vert}$ and
assume that $v_{n}\rightharpoonup v_{0}$ in $H^{1}(\Omega)$, $v_{n}\rightarrow
v_{0}$ in $L^{t}(\Omega)$ for $1\leq t<2^{\ast}$ and in $L^{2}(\partial
\Omega)$, and $v_{n}\rightarrow v_{0}$ a.e. in $\Omega$, for some $v_{0}\in
H^{1}(\Omega)$. Then
\[
\int_{\Omega}|\nabla v_{0}|^{2}-\alpha\int_{\partial\Omega}v_{0}^{2}%
\leq\liminf\left(  \int_{\Omega}|\nabla v_{n}|^{2}-\alpha\int_{\partial\Omega
}v_{n}^{2}\right)  \leq0
\]
and $\int_{\Omega}a(x)|v_{n}|^{q+1}=\Vert u_{n}\Vert^{-(q+1)}\rightarrow0$.
Hence $\int_{\Omega}a(x)|v_{0}|^{q+1}=0$. Moreover, $v_{0}\not \equiv 0$ since
otherwise, from the above inequality, we would have $v_{n}\rightarrow0$ in
$H^{1}(\Omega)$, which is impossible. Thus we have $\alpha\geq\frac
{\int_{\Omega}|\nabla v_{0}|^{2}}{\int_{\partial\Omega}v_{0}^{2}}$, which
contradicts $\alpha<\tilde{\alpha}$. Therefore $\mu(\alpha)$ is finite, and
repeating the above argument we can show that it is achieved by some
nonnegative $u$. By the Lagrange multipliers rule, we find that $u$ satisfies
$-\Delta u=\mu(\alpha)a(x)u^{q}$ in $\Omega$ and $\partial_{\nu}u=\alpha u$ on
$\partial\Omega$. Note that since $\alpha<\tilde{\alpha}$ we have $\mu
(\alpha)>0$. We set $u_{\alpha}:=\mu(\alpha)^{-\frac{1}{1-q}}u$ to get a
%nonnegative
solution of $(P_{\alpha})$ such that $\int_{\Omega}a(x)u_{\alpha}^{q+1}>0$, so
that $u_{\alpha}\not \equiv 0$ in $\Omega_{a}^{+}$. Now, if $q\in
\mathcal{A}_{D}$ then, from Proposition \ref{cor:Aal:ordered} it follows that
$q\in\mathcal{A}_{\alpha}$ for every $\alpha<0$, so that $u_{\alpha}%
\in\mathcal{P}^{\circ}$. Since $(P_{\alpha})$ has at most one solution in
$\mathcal{P}^{\circ}$ for each $\alpha<0$ (see Theorem \ref{tin}(iii)), we
infer that $u_{\alpha}$ is the unique nontrivial solution of $(P_{\alpha})$.

Finally, assume $(A.0)$ and $q\in\mathcal{A}_{N}$. Then, for $0<\alpha
<\tilde{\alpha}$ we have that $u_{\alpha}$ is a supersolution of $\left(
P_{0}\right)  $. Thus, since it is easy to provide small nontrivial
subsolutions of $\left(  P_{0}\right)  $ (see e.g. the construction in
\cite[Lemma 2.3(ii)]{BPT1}), recalling Theorem \ref{tin}(v) we deduce that
$u_{\alpha}\geq u_{N}$ on $\overline{\Omega}$, and we get the desired
conclusion. \qed\newline

Next, for $\alpha\leq0$ and $u=u_{\alpha}$, we consider the eigenvalue
problem
\[%
\begin{cases}
-\Delta\phi=a(x)qu^{q-1}\phi+\gamma\phi & \mbox{in $\Omega$},\\
\partial_{\nu}\phi=\alpha\phi & \mbox{on $\partial \Omega$},
\end{cases}
\]
where $\gamma=\gamma(\alpha,u)$ is an eigenvalue parameter. It is well known
that this problem has a
%\marginpar{modified}
smallest eigenvalue $\gamma_{1}=\gamma_{1}(\alpha,u)$, which is simple and
possesses an eigenfunction $\phi_{1}\in\mathcal{P}^{\circ}$.

\begin{lem}
[Non-degeneracy]\label{prop:nondege} Whenever $u_{\alpha}$ exists for
$\alpha\leq0$, we have $\gamma_{1}(\alpha,u_{\alpha})>0$.
\end{lem}

\noindent\textit{Proof.} By a direct computation and using Green's formula we
infer that
\begin{align*}
\gamma_{1}\int_{\Omega}u^{q}\phi_{1}  &  =\int_{\Omega}\left(  -\Delta\phi
_{1}\cdot u^{q}+\Delta u\cdot qu^{q-1}\phi_{1}\right) \\
&  =\int_{\Omega}\nabla\phi_{1}\cdot\nabla(u^{q})-\nabla u\cdot\nabla\left(
qu^{q-1}\phi_{1}\right)  +\int_{\partial\Omega}qu^{q-1}\phi_{1}\partial_{\nu
}u-u^{q}\partial_{\nu}\phi_{1}\\
&  =q(1-q)\int_{\Omega}u^{q-2}|\nabla u|^{2}\phi_{1}-\alpha(1-q)\int
_{\partial\Omega}u^{q}\phi_{1},
\end{align*}
and the conclusion follows. \qed

\begin{prop}
[Existence of an increasing curve]\label{pc} Assume $(A.1)$ and $q\in
\mathcal{A}_{D}$. Then $\alpha\mapsto u_{\alpha}$ is $C^{\infty}$ from
$(-\infty,0)$ into $W^{2,r}(\Omega)$ and $u_{\alpha}<u_{\beta}$ on
$\overline{\Omega}$ for $\alpha<\beta<0$. Moreover $u_{\alpha}\rightarrow
u_{D}$ in $H^{1}(\Omega)$ as $\alpha\rightarrow-\infty$.
\end{prop}

\noindent\textit{Proof.} Based on Lemma \ref{prop:nondege}, we show that
$\alpha\mapsto u_{\alpha}$ is $C^{\infty}$ from $(-\infty,0)$ into
$W^{2,r}(\Omega)$. Let $\delta>0$ and $B_{0}$ be a small open ball in
$W^{2,r}(\Omega)$ with center $u_{\alpha}$, so that $B_{0}\subset
\mathcal{P}^{\circ}$. Set
\begin{align*}
\mathcal{F}:  &  (\alpha-\delta,\alpha+\delta)\times B_{0}\longrightarrow
L^{r}(\Omega)\times W^{1-\frac{1}{r},r}(\partial\Omega)\\
&  \qquad(\beta,u)\longmapsto\left(  -\Delta u-a(x)u^{q},\;\partial_{\nu}
u-\beta u\right)  .
\end{align*}
We see that $\mathcal{F}(\alpha,u_{\alpha})=0$, and the Fr\'{e}chet derivative
$\mathcal{F}_{u}(\beta,u):W^{2,r}(\Omega)\rightarrow L^{r}(\Omega)\times
W^{1-\frac{1}{r},r}(\partial\Omega)$ is given by $\mathcal{F}_{u}%
(\beta,u)\varphi=\left(  -\Delta\varphi-a(x)qu^{q-1} \varphi, \partial_{\nu
}\varphi-\beta\varphi\right)  $.
%Since $u \in\mathcal{P}^{\circ}$, we have $aqu^{q-1} \varphi\in C(\overline{\Omega}) \subset L^{r}(\Omega)$.
From Lemma \ref{prop:nondege} we infer that $\mathcal{F}_{u}(\alpha,u_{\alpha
})$ is a homeomorphism, using the index theory for Fredholm operators, and
thus, the desired assertion follows by the implicit function theorem.
%Thanks to Lemma \ref{prop:nondege}, we can apply the implicit function theorem at $(\alpha_{1},u_{\alpha_{1}})$ for every $\alpha_{1}<0$, and thus we obtain locally a $C^{\infty}$ curve $(\alpha,u_{\alpha})$ near $(\alpha_{1},u_{\alpha_{1}})$.

We may then differentiate $(P_{\alpha})$ with respect to $\alpha$ to obtain
\[%
\begin{cases}
-\Delta u_{\alpha}^{\prime}=qa(x) u_{\alpha}^{q-1}u_{\alpha}^{\prime} &
\mbox{in $\Omega$},\\
\partial_{\nu} u_{\alpha}^{\prime}=u_{\alpha}+\alpha u_{\alpha}^{\prime} &
\mbox{on $\partial \Omega$}.
\end{cases}
\]
Set $\mathcal{L}:=-\Delta-qa(x) u_{\alpha}^{q-1}$ and $\mathcal{B}%
:=\frac{\partial}{\partial\nu}-\alpha$. It follows that
\[%
\begin{cases}
\mathcal{L}u_{\alpha}^{\prime}=0 & \mbox{in $\Omega$},\\
\mathcal{B}u_{\alpha}^{\prime}=u_{\alpha}>0 & \mbox{on $\partial \Omega$}.
\end{cases}
\]
Lemma \ref{prop:nondege} enables us to apply \cite[Theorem 7.10]{LGbook13} to
deduce that $u_{\alpha}^{\prime}\in\mathcal{P}^{\circ}$ for every $\alpha<0$,
which shows that $u_{\alpha}$ is increasing with respect to $\alpha$.

Let now $\alpha_{n}\rightarrow-\infty$ and $u_{n}:=u_{\alpha_{n}}$. We may
assume that $\alpha_{n}$ is decreasing, and so is $u_{n}$. Thus $\Vert
u_{n}\Vert_{C(\overline{\Omega})}$ is clearly bounded, and since $u_{n}$ is a
solution of $(P_{\alpha_{n}})$, we deduce that $\Vert u_{n}\Vert$ is bounded.
Hence, up to a subsequence, $u_{n}\rightharpoonup u_{\infty}$ in $H^{1}%
(\Omega)$, $u_{n}\rightarrow u_{\infty}$ in $L^{t}(\Omega)$ for $1\leq
t<2^{\ast}$, and in $L^{2}(\partial\Omega)$, and $u_{n}\rightarrow u_{\infty}$
a.e. in $\Omega$, for some $u_{\infty}\in H^{1}(\Omega)$. In particular,
$u_{\infty}$ is nonnegative.
%Moreover,
%\[
%\int_{\Omega}|\nabla u_{n}|^{2} = \int_{\Omega}a(x)u_{n}^{q+1} + \alpha_{n} \int_{\partial\Omega} u_{n}^{2}.
%\]
Since $u_{n}$ is a solution of $(P_{\alpha_{n}})$, we obtain
%\marginpar{argument simpified}%
\[
-\alpha_{n}\int_{\partial\Omega}u_{n}=\int_{\Omega}a(x)u_{n}^{q}\leq C\Vert
u_{n}\Vert_{C(\overline{\Omega})}^{q}\leq C^{\prime}.
\]
As $-\alpha_{n}\rightarrow\infty$, it follows that $\int_{\partial\Omega}%
u_{n}\rightarrow\int_{\partial\Omega}u_{\infty}=0$, so that $u_{\infty}=0$ on
$\partial\Omega$, implying $u_{\infty}\in H_{0}^{1}(\Omega)$. Using the
different convergences of $u_{n}$ towards $u_{0}$ and standard arguments, we
find that $u_{n}\rightarrow u_{\infty}$ in $H^{1}(\Omega)$. From the weak
formulation of $(P_{\alpha})$ we deduce that $u_{\infty}$ is a weak solution
of $(P_{-\infty})$. Finally, note from \eqref{mua} that $\mu(\alpha)\leq
\int_{\Omega}|\nabla v|^{2}$ for any $v\in H_{0}^{1}(\Omega)$ such that
$\int_{\Omega}a(x)|v|^{q+1}=1$. Hence $\mu(\alpha)\leq C$ for some constant
$C>0$ and any $\alpha$. It follows that
%\marginpar{space reduced}
$\int_{\Omega}a(x)u_{n}^{q+1}=\mu(\alpha_{n})^{-\frac{q+1}{1-q}}\geq
C^{-\frac{q+1}{1-q}}>0$ for every $n$, which implies that $u_{\infty}$ is
nontrivial. Since $q\in\mathcal{A}_{D}$ we have $u_{\infty}=u_{D}$, as
desired. \qed\newline

\begin{prop}
[Asymptotic behavior as $\alpha\to0^-$]\label{last}Assume $(A.1)$ and
$q\in\mathcal{A}_{D}$.

\begin{enumerate}
\item If $\int_{\Omega}a\geq0$ then $\displaystyle\min_{x\in\overline{\Omega}%
}u_{\alpha}(x)\rightarrow\infty$ as $\alpha\rightarrow0^{-}$, and $(P_{\alpha
})$ has no solution $u$ such that $u\not \equiv 0$ in $\Omega_{+}^{a}$ for
$\alpha\geq0$ (in particular it has no nontrivial solution for $0\leq
\alpha<\alpha_{p}$).

\item If $\int_{\Omega}a <0$ then the curve $\alpha\mapsto u_{\alpha}$ can be
extended to $(-\infty,\overline{\alpha})$, for some $\overline{\alpha}>0$, so
that $u_{0}=u_{N}$ and $u_{\alpha}\in\mathcal{P}^{\circ}$ is a solution of
$(P_{\alpha})$ for $\alpha\in(0,\overline{\alpha})$. Moreover, $\alpha\mapsto
u_{\alpha}$ is increasing in $(-\infty,\overline{\alpha})$, and unique in the
following sense: if $u_{n}$ is a solution of $(P_{\alpha_{n}})$ such that
$\alpha_{n}\rightarrow0^{+}$ and $u_{n} \rightarrow u_{N}$ in $C^{1}%
(\overline{\Omega})$, then, for $n$ large enough, $u_{n}=u_{\alpha}$ for some
$\alpha\in(0,\overline{\alpha})$.
\end{enumerate}
\end{prop}

\noindent\textit{Proof.}

\begin{enumerate}
\item First we prove $\Vert u_{\alpha}\Vert_{C(\overline{\Omega})}%
\rightarrow\infty$ as $\alpha\rightarrow0^{-}$. Since $u_{\alpha}%
\in\mathcal{P}^{\circ}$ is a solution of $(P_{\alpha})$, it suffices to show
$\Vert u_{\alpha}\Vert\rightarrow\infty$ as $\alpha\rightarrow0^{-}$. Assume
by contradiction that for some sequence $\alpha_{n}\rightarrow0^{-}$, $\Vert
u_{\alpha_{n}}\Vert$ is bounded. By elliptic regularity, it follows that, up
to a subsequence, $u_{\alpha_{n}}\rightarrow u_{\ast}$ in $C^{1}%
(\overline{\Omega})$ for some $u_{\ast}$. By definition, we deduce that
$u_{\ast}$ is a solution of $(P_{0})$. Moreover, $u_{\ast}\in\mathcal{P}%
^{\circ}$ by the monotonicity of $u_{\alpha_{n}}$, i.e. $u_{\ast}=u_{N}$.
%
%Furthermore, we can fix a ball $B\subset\Omega_{a}^{+}$, so that, since $u_{\alpha_{n}}>c$ in $B$ for some $c>0$  \cite[Lemma 2.2]{KRQU16}, we have that $u_{\ast}\not \equiv 0$. Lastly, since $q\in\mathcal{A}_{N}$ (by Proposition  \ref{cor:Aal:ordered}), we infer that $u_{\ast}\in\mathcal{P}^{\circ}$ is a positive solution of $(P_{0})$, i.e. $u_{\ast}=u_{N}$.
However, this contradicts \cite[Corollary 2.1]{BPT2} (which clearly holds in
our setting), as desired.

Now, by monotonicity it suffices to show the existence of a sequence
$\alpha_{n}\rightarrow0^{-}$ such that $\min_{\overline{\Omega}}u_{\alpha_{n}%
}\rightarrow\infty$. Let $\alpha_{n}\rightarrow0^{-}$. Set $u_{n}%
:=u_{\alpha_{n}}$ and $v_{n}:=u_{n}/\Vert u_{n}\Vert_{C(\overline{\Omega})}$.
Then it follows that
\[
\int_{\Omega}|\nabla v_{n}|^{2}\leq\left(  \int_{\Omega}a(x)v_{n}%
^{q+1}\right)  \Vert u_{n}\Vert_{C(\overline{\Omega})}^{q-1}\rightarrow0.
\]
We deduce that, up to a subsequence, $v_{n}\rightarrow c_{\ast}$ in
$H^{1}(\Omega)$, where $c_{\ast}$ is a nonnegative constant. Since $v_{n}$
satisfies
\[%
\begin{cases}
-\Delta v_{n}=a(x)\Vert u_{n}\Vert_{C(\overline{\Omega})}^{q-1}v_{n}^{q} &
\mbox{in $\Omega$},\\
\partial_{\nu} v_{n}=\alpha_{n}v_{n} & \mbox{on $\partial \Omega$},
\end{cases}
\]
we find that $\Vert v_{n}\Vert_{C^{\theta}(\overline{\Omega})}$ is bounded for
$\theta\in(0,1)$ by elliptic regularity and a bootstrap argument \cite[Theorem
2.2]{Ro05}. By a compactness argument, we infer that, up to a subsequence,
$v_{n}\rightarrow c_{\ast}$ in $C(\overline{\Omega})$ and $c_{\ast}>0$, from
which our desired conclusion follows.

Finally, if $u$ is a nontrivial solution of $(P_{\alpha})$ such that
$u\not \equiv 0$ in $\Omega_{+}^{a}$ and $\alpha\geq0$ then $u$ is a
supersolution of $(P_{0})$. Hence $(P_{0})$ has a nontrivial solution $u_{0}$,
and since $q\in\mathcal{A}_{D}$, we have $u_{0}\in\mathcal{P}^{\circ}$.
Reasoning as in \cite[Lemma 2.1]{BPT2} we infer that $(A.0)$ holds, which
contradicts our assumption.

\item From $(A.0)$ and $q\in\mathcal{A}_{N}$ (by Proposition
\ref{cor:Aal:ordered}), we know that $u_{N}\in\mathcal{P}^{\circ}$ is the
unique nontrivial solution of $(P_{0})$. By Lemma \ref{prop:nondege} we have
$\gamma_{1}(0,u_{N})>0$. Arguing as in the proof of Proposition \ref{pc}, the
implicit function theorem allows us to find some $\overline{\alpha}>0$ and an
increasing $C^{\infty}$ curve $(\alpha,u_{\alpha})$ with $u_{\alpha}%
\in\mathcal{P}^{\circ}$ solutions of $(P_{\alpha})$, parametrized by
$\alpha\in(-\overline{\alpha},\overline{\alpha})$. Lastly, let $u_{n}$ be a
nontrivial solution of $(P_{\alpha_{n}})$ such that $\alpha_{n}\rightarrow
0^{+}$ and $u_{n}\rightarrow u_{0}$ in $C^{1}(\overline{\Omega})$. So, the
Lebesgue dominated convergence theorem shows that
\[%
\begin{cases}
-\Delta(u_{n}-u_{0})=a(x)(u_{n}^{q}-u_{0}^{q})\rightarrow0 &
\mbox{in  $L^r(\Omega)$},\\
\partial_{\nu}(u_{n}-u_{0})+(u_{n}-u_{0})=\alpha_{n}u_{n}+(u_{n}%
-u_{0})\rightarrow0 & \mbox{in  $C^1(\partial \Omega)$}.
\end{cases}
\]
We deduce then, by elliptic regularity, that $u_{n}\rightarrow u_{0}$ in
$W^{2,r}(\Omega)$. Combining the existence result with an application of the
implicit function theorem provides the desired assertion. \qed\newline
\end{enumerate}

%\end{enumerate}
%%
%letting $\alpha_{n} \rightarrow-\infty$, we have, up to a subsequence, $u_{\alpha_{n}}\rightarrow u_{\infty}$ in $H^{1}(\Omega)$ for some nontrivial solution $u_{\infty}$ of $(P_{-\infty})$ in $H_{0}^{1}(\Omega)\cap  W^{2,r}(\Omega)$ with $r>N$. By using $q\in\mathcal{A}_{D}$, it follows that $u_{\infty}\in\mathcal{P}^{\circ}$. Thanks to the uniqueness of a solution in
%$\mathcal{P}^{\circ}$ for $(P_{-\infty})$, we conclude that $u_{\alpha
%}\rightarrow u_{D}$ in $H^{1}(\Omega)$. The proof is now complete. \qed

%Section

%\section{Bifurcation analysis of $(R_\alpha)$}

\section{Qualitative analysis and exact multiplicity}

\label{secpos}

%\subsection{Auxiliary results for $(P_{\alpha})$}

%\label{subsec:aP}
%\subsection{Local bifurcation analysis for $(P_{\alpha})$ via the ImplicitFunction Theorem}

%\label{subsec:Q}

%\subsection{The rescaled problem $(R_{\alpha})$}

%\label{subsec:R}

\begin{figure}[tbh]
\begin{center}
\includegraphics[scale=0.18]{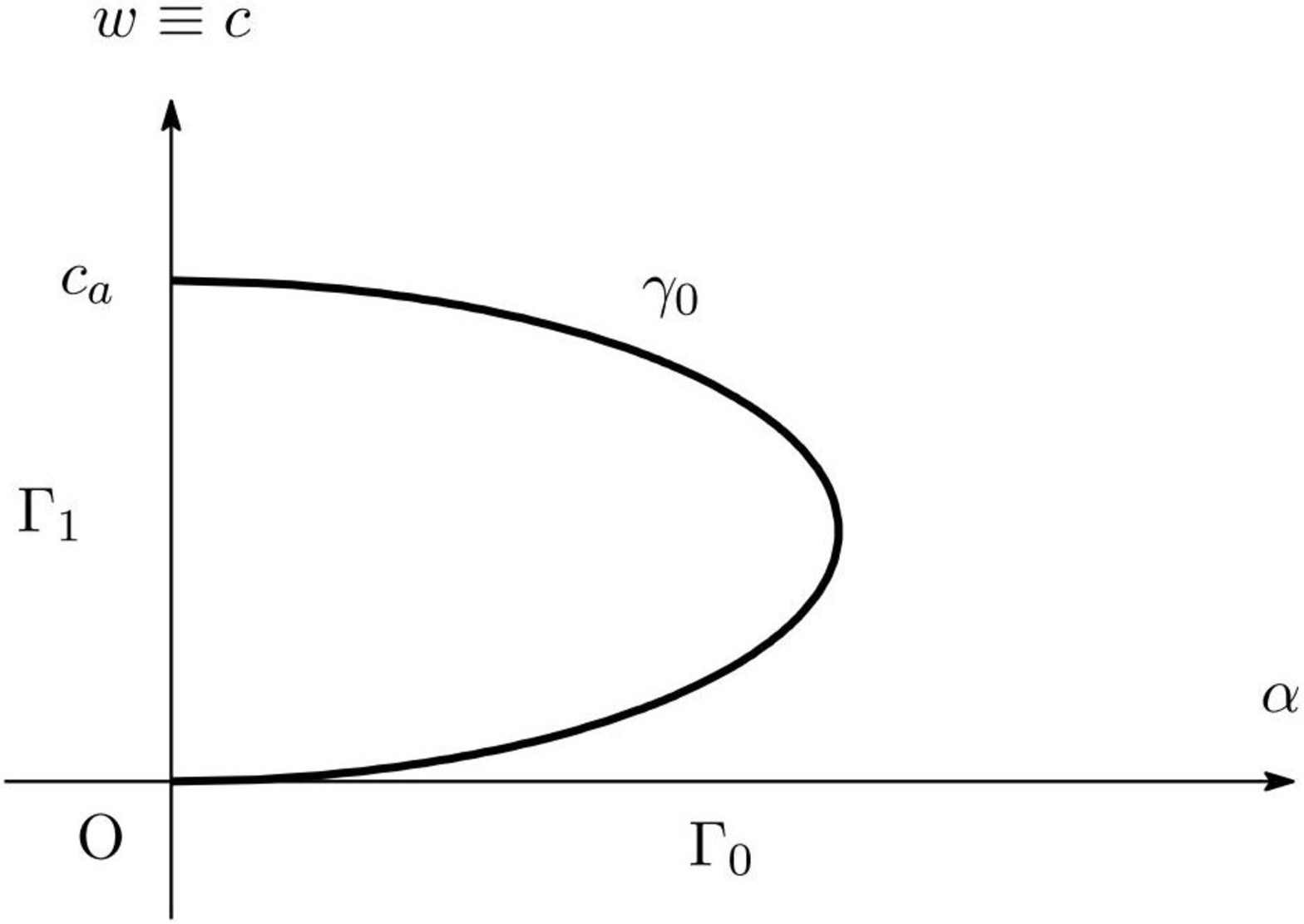}
\end{center}
\caption{The bounded component $\gamma_{0}$ of nontrivial solutions of
$(R_{\alpha})$ when $\int_{\Omega}a<0$.}%
\label{fig18_0105b}%
\end{figure}

In this section we prove an exact multiplicity result for $(R_{\alpha})$.
Furthermore, we establish some preliminary results to prove Theorem
\ref{maint} below, which states the existence of a subcontinuum $\gamma
_{0}=\{(\alpha,w)\in\lbrack0,\infty)\times C^{1}(\overline{\Omega})\}$ of
solutions of $(R_{\alpha})$ such that
\begin{equation}
\gamma_{0}\cap(\Gamma_{0}\cup\Gamma_{1})=\left\{  (0,0),(0,c_{a})\right\}
\label{gcapGab}%
\end{equation}
(recall that $\Gamma_{0}$ and $\Gamma_{1}$ are the solution lines of $\left(
R_{\alpha}\right)  $ given by (\ref{G0:00}), see Figure \ref{fig18_0105b}).
%with $\alpha= \alpha$ and $c_{a}$ replaced by $c_a$
%{fig18_0106a}
We shall use this result to prove Theorem \ref{thm:main}(iii).
%	% \begin{figure}[H]

%	% \begin{figure}[H]
%\begin{figure}[tbh]
%\begin{center}
%\includegraphics[scale=0.2]{fig18_0106a.eps}
%\end{center}
%\caption{A minimal possible bounded component $\gamma_{0}$.}%
%\label{fig18_0106a}%
%\end{figure}

%\subsection{A priori upper and lower bounds and preliminary bifurcation analysis}

%\label{sec:bd}

%\subsubsection{Bifurcation from $\Gamma_1$}

First we show the existence of an \textit{a priori} lower bound in
$C(\overline{\Omega})$ for positive solutions of $(P_{\alpha})$ with
$\alpha\in\mathbb{R}$, which shows that such solutions do not bifurcate from
zero at any $\alpha\in\mathbb{R}$:

\begin{lem}
[A priori lower bound]\label{prop:glb} There exists $C>0$ such that
$\|u\|_{C(\overline{\Omega})} \geq C$ for every positive supersolution of
$(P_{\alpha})$ and every $\alpha\in\mathbb{R}$. In particular, given $\beta>0$
there exists $C_{\beta}>0$
%\marginpar{modified}
such that $\|w\|_{C(\overline{\Omega})} \geq C_{\beta}$ for every positive
supersolution of $(R_{\alpha})$ with $\alpha\geq\beta$.
\end{lem}

\noindent\textit{Proof.} The first inequality is a direct consequence of
\cite[Lemma 2.2]{KRQU16}, and by the change of variables \eqref{changeofv}, we
see that it implies the second one.\qed\newline

%-----------------------------------------

%\subsubsection{Bifurcation from infinity and a priori upper bounds}

Second we discuss bifurcation from infinity at $\alpha\geq0$. The following
result asserts that $\alpha=0$ is the only point where solutions of
$(R_{\alpha})$ bifurcate from infinity, and such solutions are given precisely
by $\Gamma_{1}$.

%-------------------------------------

%------------------------------------------------------------

\begin{prop}
[Bifurcation from infinity and a priori upper bounds]\label{cor:bounds} Given
$\overline{\alpha}>0$, there exists $C_{\overline{\alpha}}>0$ such that $\Vert
w\Vert_{C^{1}(\overline{\Omega})}\leq C_{\overline{\alpha}}$ for all
%nonnegative
solutions $w$ of $(R_{\alpha})$ with $0<\alpha<\overline{\alpha}$.
\end{prop}

%-------------------------------------

\noindent\textit{Proof.} Assume by contradiction that there exist $(\alpha
_{n},w_{n})$ such that $w_{n}$ is a solution of $(R_{\alpha_{n}})$,
$\alpha_{n}\rightarrow\alpha_{0}\geq0$, and $\Vert w_{n}\Vert_{C^{1}%
(\overline{\Omega})}\rightarrow\infty$. By elliptic regularity, it follows
that $\Vert w_{n}\Vert\rightarrow\infty$. If we set $\psi_{n}:=w_{n}/\Vert
w_{n}\Vert$ then we may assume that for some $\psi\in H^{1}(\Omega)$ and
$1\leq t<2^{\ast}$,%
\begin{equation}
\left\{
\begin{array}
[c]{l}%
\psi_{n}\rightharpoonup\psi\quad\mbox{in}\quad H^{1}(\Omega),\qquad\psi
_{n}\rightarrow\psi\quad\mbox{in}\quad L^{t}(\Omega),\\
\psi_{n}\rightarrow\psi\quad\mbox{in}\quad L^{2}(\partial\Omega),\qquad
\psi_{n}\rightarrow\psi\quad\mbox{a.e. in $\Omega$.}
\end{array}
\right.  \label{psiconverg}%
\end{equation}
Since $w_{n}$ is a weak solution of $(R_{\alpha_{n}})$, we see that
\[
\int_{\Omega}\nabla w_{n}\nabla\varphi=\alpha_{n}\int_{\Omega}a\left(
x\right)  w_{n}^{\;q}\varphi+\alpha_{n}\int_{\partial\Omega}w_{n}\varphi
,\quad\forall\varphi\in H^{1}(\Omega).
\]
Dividing it by $\Vert w_{n}\Vert$, it follows that
\begin{equation}
\int_{\Omega}\nabla\psi_{n}\nabla\varphi=\alpha_{n}\left(  \int_{\Omega
}a(x)\psi_{n}^{\;q}\varphi\right)  \Vert w_{n}\Vert^{q-1}+\alpha_{n}%
\int_{\partial\Omega}\psi_{n}\varphi, \label{un2.4}%
\end{equation}
so that $\int_{\Omega}\nabla\psi\nabla\varphi=\alpha_{0}\int_{\partial\Omega
}\psi\varphi$ for all $\varphi\in H^{1}(\Omega)$. Hence $\psi\geq0$ solves the
problem
\begin{equation}%
\begin{cases}
\Delta\psi=0 & \mbox{in $\Omega$},\\
\partial_{\nu}\psi=\alpha_{0}\psi & \mbox{on $\partial \Omega$},
\end{cases}
\label{epro0b}%
\end{equation}
Taking $\varphi\equiv1$ in \eqref{un2.4} we find that $\Vert w_{n}\Vert
^{q-1}\int_{\Omega}a(x)\psi_{n}^{\;q}+\int_{\partial\Omega}\psi_{n}=0$.
Passing to the limit, we have $\int_{\partial\Omega}\psi=0$, i.e. $\psi
\equiv0$ on $\partial\Omega$, so that $\psi\equiv0$ from \eqref{epro0b}. Since
$\Vert\psi_{n}\Vert=1$, we deduce that $\psi_{n}\rightharpoonup0$
%\mpfns{\bf two modifications}
but $\psi_{n}\not \rightarrow 0$ in $H^{1}(\Omega)$.

Finally, taking $\varphi=\psi_{n}$ in \eqref{un2.4} we find that
\[
0=\int_{\Omega}|\nabla\psi|^{2}\leq\displaystyle\liminf_{n\rightarrow\infty
}\int_{\Omega}|\nabla\psi_{n}|^{2}=\alpha_{0}\int_{\partial\Omega}\psi^{2}=0,
\]
and thus that $\psi_{n}\rightarrow0$ in $H^{1}(\Omega)$, a contradiction.
\qed\newline

%------------------------------------------------------------

%\subsubsection{A priori bounds on $\alpha$}

\label{subsec:ubal}

Using Proposition \ref{cor:bounds} we show that (under the conditions of
Theorem \ref{thm:main}) the existence range for nontrivial solutions of
$(P_{\alpha})$ is an interval. We set
\[
\alpha_{s}:=\sup\{\alpha>0:\mbox{$(P_\alpha)$ has a nontrivial solution}\}.
\]
Note that this definition is equivalent to \eqref{as} if $(A.2)$ holds
%\marginpar{\textit{modified}}
and $q\in\mathcal{A}_{N}$, in view of Lemma \ref{lpos} and Proposition
\ref{cor:Aal:ordered}.

\begin{prop}
\label{prop:exnon} Assume $(A.1)$ and $(A.2)$. If $(P_{\alpha})$ has a
nontrivial solution for some $\alpha>0$, then
%$\alpha_{s}\in(0,\infty]$, and
$(P_{\alpha})$ has at least one nontrivial solution for $0\leq\alpha\leq
\alpha_{s}$ ($0\leq\alpha<\infty$ if $\alpha_{s}=\infty$).
%\marginpar{removed, since holds by definition of $\alpha_{s}$}.

\end{prop}

\noindent\textit{Proof.} We may assume that $\alpha_{s}<\infty$. Then
$(P_{\alpha_{s}})$ has a nontrivial solution $u_{s}$ by elliptic regularity,
using Lemma \ref{prop:glb} and Proposition \ref{cor:bounds}. In this case,
$u_{s}$ is a supersolution of $(P_{\alpha})$ for every $\alpha<\alpha_{s}$ and
$u_{s} \not \equiv 0$ in $\Omega_{+}^{a}$ by Lemma \ref{lpos}. We can now
deduce that $(P_{\alpha})$ has at least one nontrivial solution for each
$\alpha\in\lbrack0,\alpha_{s})$ by constructing a suitable small subsolution
(see the proof of Proposition \ref{cor:Aal:ordered}),
%\marginpar{modified},
as desired. \qed\newline

Third we establish an \textit{a priori} bound on $\alpha>0$ for the existence
of solutions in $\mathcal{P}^{\circ}$ of $(P_{\alpha})$ and $(R_{\alpha})$.

\begin{prop}
[A priori bounds on $\alpha$ for $q\in\mathcal{A}_{N}$]\label{alpha}Assume
$\left(  A.0\right)  $, $\left(  A.1\right)  $ and $q\in\mathcal{A}_{N}$. If
$\left(  P_{\alpha}\right)  $ or $\left(  R_{\alpha}\right)  $ has a
supersolution in $\mathcal{P}^{\circ}$ with $\alpha>0$, then $\alpha
<\frac{-\int_{\Omega}a}{\int_{\partial\Omega}u_{N}^{1-q}}$, where $u_{N}%
\in\mathcal{P}^{\circ}$ is the unique nontrivial solution of $\left(
P_{0}\right)  $.
\end{prop}

\noindent\textit{Proof. }Taking into account the change of variables
\eqref{changeofv}, we consider without loss of generality the problem $\left(
P_{\alpha}\right)  $. Suppose $\left(  P_{\alpha}\right)  $ has a
supersolution $u_{\alpha}\in\mathcal{P}^{\circ}$. Then $u_{\alpha}$ is a
supersolution of $\left(  P_{0}\right)  $. Using a suitably small first
eigenfunction (under homogeneous Dirichlet boundary condition) with respect to
the weight $a$ in some smooth subdomain of $\Omega_{+}^{a}$ and extending it
by zero to $\Omega$, we obtain a nontrivial weak subsolution of $\left(
P_{0}\right)  $ smaller than $u_{\alpha}$. Hence, we get a nontrivial solution
$v$ of $\left(  P_{0}\right)  $, with $v\leq u_{\alpha}$ in $\Omega$. Now,
since $q\in\mathcal{A}_{N}$, from Theorem \ref{tin}(v) we deduce that
$v=u_{N}\in\mathcal{P}^{\circ}$.\newline On the other hand, taking $u_{\alpha
}^{-q}$ as a test function in the weak form of $\left(  P_{\alpha}\right)  $
we get that
\[
-q\int_{\Omega}\frac{\left\vert \nabla u_{\alpha}\right\vert ^{2}}{u_{\alpha
}^{q+1}}\geq\int_{\Omega}a+\alpha\int_{\partial\Omega}u_{\alpha}^{1-q}.
\]
Therefore,%
\[
\alpha\int_{\partial\Omega}u_{N}^{1-q}\leq\alpha\int_{\partial\Omega}%
u_{\alpha}^{1-q}<-\int_{\Omega}a
\]
and the conclusion follows. \qed\newline

When $q\in\left[  0,1\right)  $ we can still provide an \textit{a priori}
bound similar to the previous one.
%
%We remark that $(H)$ is weaker than $(A.3)$ (see Figures \ref{figA23}(ii) and
%\ref{figsubdD02}).
%To stress the dependence on $a$, we write $(P_{\alpha,a})$ instead of $(P_\alpha)$.
%\begin{figure}[tbh]
%\centerline{
%\includegraphics[scale=0.14]{fig18_1202H.eps}
%\hskip0.35cm
%\includegraphics[scale=0.14] {fig18_1202A3.eps}
%} \centerline{(i) \hskip6.0cm (ii) }
%\caption{(i) Condition $(H)$; (ii) Condition $(A.3)$.}
%\label{figsubdD02}
%\end{figure}
%	% \begin{figure}[H]
%\begin{figure}[tbh]
%\begin{center}
%\includegraphics[scale=0.18]{figsubdD02.eps}
%\end{center}
%\caption{The subdomain $D$ and ball $B$.} \label{figsubdD02}
%\end{figure}
Before stating this result, we need to establish the uniqueness of positive
solutions for the following concave mixed problem:
\begin{equation}%
\begin{cases}
-\Delta u=a(x)f(u) & \mbox{in $D$},\\
u=0 & \mbox{on $\Gamma_{\Omega}$},\\
\partial_{\nu}u=0 & \mbox{on $\Gamma_{\partial \Omega}$},
\end{cases}
\label{cmpD}%
\end{equation}
where $f:[0,\infty)\rightarrow\mathbb{R}$ is continuous, and $f(s)/s$ is
decreasing for $s>0$. Recall that $D$, $\Gamma_{\Omega}$, and $\Gamma
_{\partial\Omega}$ are given by $(A.3)$.

\begin{lem}
\label{lem:ump} Assume $(A.3)$. Then \eqref{cmpD} has at most one positive solution.
\end{lem}

\noindent\textit{Proof.} Let $u_{1},u_{2}$ be positive solutions of
\eqref{cmpD}. Then, for $i=1,2$ we have $u_{i}\in H_{\Gamma_{\Omega}}^{1}(D)$
and
\begin{equation}
\int_{D}\nabla u_{i}\nabla\varphi=\int_{D}a(x)f(u_{i})\varphi\quad
\mbox{ for all }\varphi\in H_{\Gamma_{\Omega}}^{1}(D), \label{u1u2}%
\end{equation}
where $H_{\Gamma_{\Omega}}^{1}(D):=\left\{  u\in H^{1}(D):u=0\mbox{ on }\Gamma
_{\Omega}\right\}  $. Arguing as in the proof of \cite[Proposition
A.1]{RQUTMNA}, we deduce that for $i\neq j$,
\[
\int_{\{u_{i}>u_{j}\}}u_{i}u_{j}a(x)\left(  \frac{f(u_{j})}{u_{j}}%
-\frac{f(u_{i})}{u_{i}}\right)  \leq0.
\]
It follows that $u_{1}=u_{2}$ in $\{x\in D:a(x)>0\}$, so $a(x)f(u_{1}%
)=a(x)f(u_{2})$ in $D$. Going back to \eqref{u1u2}, we deduce the desired
conclusion. \qed\newline

\begin{prop}
[A priori bounds on $\alpha$ for $q\in\lbrack0,1)$]\label{alpha2}Assume
$\left(  A.0\right)  $, $\left(  A.3\right)  $ and $q\in\left[  0,1\right)  $.
If $\left(  P_{\alpha}\right)  $ or $\left(  R_{\alpha}\right)  $ has a
supersolution in $\mathcal{P}^{\circ}$ with $\alpha>0$, then $\alpha
<\frac{-\int_{\Omega}a}{\int_{\Gamma_{\partial\Omega}}v^{1-q}}$, where $v$ is
the unique positive solution of
\[%
\begin{cases}
-\Delta v=a(x)v^{q} & \text{in }D,\\
v\geq0 & \text{in }D,\\
v=0 & \text{on }\Gamma_{\Omega},\\
\partial_{\nu}v=0 & \text{on }\Gamma_{\partial\Omega}.
\end{cases}
\leqno{(Q_{q,a})}
\]

\end{prop}

\noindent\textit{Proof. }As above, we may consider only $\left(  P_{\alpha
}\right)  $. We argue as in the proof of Proposition \ref{alpha}, with some
minor changes. Let us indicate them. Let $q\in\left[  0,1\right)  $ and
suppose that $\left(  P_{\alpha}\right)  $ has a supersolution $u_{\alpha}%
\in\mathcal{P}^{\circ}$ with $\alpha>0$. Then $u_{\alpha}$ is a supersolution
of $(Q_{q,a})$. On the other side, let $D$ be as in $(A.3)$. Taking a small
first Dirichlet eigenfunction associated to the weight $a$ in $D$, we have a
subsolution of $(Q_{q,a})$ smaller than $u_{\alpha}$. Thus, by the sub and
supersolutions method under mixed boundary conditions (see e.g. \cite{LS06}),
we obtain a nontrivial
%nonnegative
solution $v$ of $(Q_{q,a})$, with $v\leq u_{\alpha}$ in $D$. Moreover, by the
strong maximum principle and Hopf's Lemma, we have $v>0$ on $D\cup
\Gamma_{\partial\Omega}$, and in particular $\int_{\Gamma_{\partial\Omega}%
}v^{1-q}>0$. We also note that $v$ does not depend on $\alpha$ (it depends on
$q$, but $q$ is fixed), since $(Q_{q,a})$ admits at most one positive solution
by Lemma \ref{lem:ump}. Now we can conclude the argument as in the proof of
Proposition \ref{alpha}, with $v$ in place of $u_{N}$. \qed\newline

\begin{rem}
\label{2.6}\strut

\begin{enumerate}
\item Let us mention that using an approximation procedure as in \cite[Lemma
2.1]{BPT2} one can see that the estimates in Proposition \ref{alpha} and
\ref{alpha2} hold for positive supersolutions (not necessarily in
$\mathcal{P}^{\circ}$) of $(P_{\alpha})$ and $(R_{\alpha})$.

%\item In addition to $(A.0)$, assume $(A.3)$. Then, with the aid of thecomparison principle \cite[Proposition A.1]{RQUTMNA} for concave problems with mixed boundary conditions (i.e.\ as in the proof of Proposition \ref{prop:nbddalph} below), we can prove the existence of an \textit{a priori} upper bound on $\alpha$ (similar to the one in Proposition \ref{alpha2}(ii)) \textit{uniformly} with respect to $q\in\lbrack0,1)$. Namely, if $(P_{\alpha})$ has a solution in $\mathcal{P}^{\circ}$ for $\alpha>0$ then $\alpha \leq-\lambda_{0}\frac{\int_{\Omega}a}{\int_{\Gamma_{\partial\Omega}}\phi_{1}}$, where $\Gamma_{\partial\Omega}$, $\lambda_{0}$ and $\phi_{1}$ are from \eqref{ep} and \eqref{p:DD}.

\item Let $w(x):=\alpha^{\frac{1}{1-q}}\sin^{r}x$ and $a_{q}(x):=r(1-r\cos
^{2}x)$ for $0\leq x\leq\pi$, where $q\in(0,1)$ and $r:=\frac{2}{1-q}$. We may
easily check that $-w^{\prime\prime}=\alpha a_{q}(x)w^{q}$, $w>0$ in $(0,\pi
)$, and $w(0)=w(\pi)=w^{\prime}(0)=w^{\prime}(\pi)=w^{\prime\prime
}(0)=w^{\prime\prime}(\pi)=0$. This example shows that if $\left(  A.3\right)
$ does not hold, then $(R_{\alpha})$ may have positive solutions for all
$\alpha>0$. Furthermore, extending $w$ by zero to $\Omega:=(-\delta,\pi
+\delta)$, for some $\delta>0$, we see that $w$ is a nontrivial solution
(which is \textit{not} a positive solution) of $(R_{\alpha})$ for \textit{any}
$\alpha>0$, no matter how we extend $a_{q}$. In particular we see that
$q\not \in \mathcal{A}_{\alpha}(a_{q})$ for every $\alpha\in\lbrack
-\infty,\infty)$. This extension shows that $(R_{\alpha})$ may have a
nontrivial solution for every $\alpha>0$, regardless of the behavior of $a$
near the boundary.
\end{enumerate}
\end{rem}

From Propositions \ref{alinf} and \ref{alpha} we obtain the following bounds
on $\alpha_{s}$ (recall that $\tilde{\alpha}$ is given by \eqref{def:tilal}):

\begin{cor}
\label{cal} Assume $(A.0)$, $(A.1)$ and $q\in\mathcal{A}_{N}$. Then
%\marginpar{space reduced}
$\tilde{\alpha}\leq\alpha_{s}\leq\frac{-\int_{\Omega}a}{\int_{\partial\Omega
}u_{N}^{1-q}}$.
\end{cor}

\begin{rem}
Under $(A.0)$ one may proceed as in the proof of Proposition \ref{alinf} to
show that
\[
\tilde{\mu}(\alpha):=\inf\left\{  \int_{\Omega}|\nabla u|^{2}-\alpha
\int_{\partial\Omega}u^{2}:u\in H^{1}(\Omega),\int_{\Omega}a(x)|u|^{q+1}%
=-1\right\}
\]
is achieved and negative for $0<\alpha<\sigma$, where
\[
\sigma:=\inf\left\{  \int_{\Omega}|\nabla v|^{2}:v\in H^{1}(\Omega
),\int_{\Omega}a(x)|v|^{q+1}=0,\int_{\partial\Omega}v^{2}=1\right\}  .
\]
The minimiser associated to $\tilde{\mu}(\alpha)$ gives rise then
to a nontrivial solution of $(P_{\alpha})$ for $0<\alpha<\sigma$. Thus, under
the assumptions of Corollary \ref{cal}, we have $\alpha_{s} \geq\sigma$. Note
that $\sigma\geq\tilde{\alpha}$.
\end{rem}

By Lemma \ref{prop:glb} we know that $(0,0)$ is the only possible bifurcation
point in $\Gamma_{0}\cup\{(0,0)\}$ for nontrivial solutions of $(R_{\alpha})$.
In this case, we show that the corresponding solution of $(P_{\alpha})$
remains bounded in $C^{1}(\overline{\Omega})$ as $\alpha\rightarrow0^{+}$. More precisely:

%\subsubsection{Growth order of small nonnegative solutions}

%Finally we give a upper growth order, in terms of $\alpha$, for nontrivial solutions $w_{\alpha}$ of $(R_{\alpha})$ converging to $0$ in $C^{1}(\overline{\Omega})$ as $\alpha\to0^{+}$.
%	Before stating our result, we prove the following the lemma.
%	\begin{lem}
%	Assume $(A.1)$.
%	There exists $C>0$ such that $\| u \|_{C^1(\overline{\Omega})} \geq C$ for all nontrivial solutions of $(R_{1,0})$.
%	\end{lem}
%	\noindent
%	\textit{Proof.}
%	Assume to the contrary that we can choose nontrivial solutions $u_n$ of $(R_{1,0})$ such that $\| u_n \|_{C^1(\overline{\Omega})} \rightarrow 0$.
%	We assert that $u_n \not\equiv 0$ in $\Omega_{+}^{a}$. Indeed, we deduce %from $(A.1)$ that
%	\begin{align*}
%	0< \int_\Omega |\nabla u_n|^2 = \int_\Omega a(x)u_n^{q+1}
%	\leq \int_{\mathrm{supp}\, a^+} a^+(x) u_n^{q+1}
%	= \int_{\Omega_{+}^{a}} a^+(x) u_n^{q+1}.
%	\end{align*}
%	Using $(A.1)$ again, we may infer that $u_n \not\equiv 0$ in some connected component $\Omega'$ of $\Omega_{+}^{a}$. So, the result \cite[Lemma 2.2]{KRQU16} applies, and leads us to a contradiction. \qed

\begin{prop}
[Bifurcation from $(0,0)$]\label{prop:C1C2} Assume $(A.0)$. If $\alpha
_{n}\rightarrow0^{+}$ and $w_{n}$ are solutions of $(R_{\alpha_{n}})$ with
$w_{n}\rightarrow0$ in $C^{1}(\overline{\Omega})$, then $\{\alpha_{n}%
^{-\frac{1}{1-q}}w_{n}\}$ is bounded in $C^{1}(\overline{\Omega})$.
%\marginpar{modified}

%\item If $(A.1)$ \marginpar{\scriptsize $(A.2)$ removed}%and $(A.2)$
%hold then $\displaystyle\inf_{n}\Vert\alpha_{n}^{-\frac{1}{1-q}}w_{n}\Vert_{C^{1}(\overline{\Omega})}>0$.\end{enumerate}

\end{prop}

\noindent\textit{Proof.} Assume by contradiction that $\alpha_{n}%
\rightarrow0^{+}$ and $w_{n}$ are solutions of $(R_{\alpha_{n}})$ such that
$\Vert w_{n}\Vert_{C^{1}(\overline{\Omega})}\rightarrow0$ but $\Vert\alpha
_{n}^{-\frac{1}{1-q}}w_{n}\Vert_{C^{1}(\overline{\Omega})}\rightarrow\infty$.
%	Since $u_n$ is a nonnegative solution of $(P_{\alpha_n})$, we see
%	\begin{align*}
%	\int_\Omega |\nabla u_n|^2 = \alpha_n \int_\Omega a(x)u_n^{q+1} + \alpha_n \int_{\partial \Omega} b(x)u_n^2.
%	\end{align*}
Then $u_{n}:=\alpha_{n}^{-\frac{1}{1-q}}w_{n}$ solves $(P_{\alpha_{n}})$, so
that
%\marginpar{space reduced}%
\begin{equation}
\int_{\Omega}|\nabla u_{n}|^{2}=\int_{\Omega}a(x)u_{n}^{\;q+1}+\alpha_{n}%
\int_{\partial\Omega}u_{n}^{\;2}. \label{goss}%
\end{equation}
Since $\Vert u_{n}\Vert_{C^{1}(\overline{\Omega})}\rightarrow\infty$, an
elliptic regularity argument enables us to infer that $\Vert u_{n}%
\Vert\rightarrow\infty$. Setting $\psi_{n}:=u_{n}/\Vert u_{n}\Vert$, we may
assume that $\psi_{n}$ satisfies \eqref{psiconverg}.
%for some $\psi_{\infty}\in H^{1}(\Omega)$,
%\begin{align*}
%&  \psi_{n}\rightharpoonup\psi_{\infty} \quad\mbox{in}\quad H^{1}(\Omega),\\
%&  \psi_{n}\rightarrow\psi_{\infty}\quad\mbox{in}\quad L^{t}(\Omega)\quad\mbox{for}\quad1\leq t<2^{\ast},\\
%&  \psi_{n}\rightarrow\psi_{\infty}\quad\mbox{in}\quad L^{2}(\partial\Omega),\\
%&  \psi_{n}\rightarrow\psi_{\infty}\quad \mbox{a.e. in $\Omega$, in particular,}\ \psi_{\infty} \geq0.
%\end{align*}
Moreover, dividing by $\Vert u_{n}\Vert^{2}$, it follows from \eqref{goss}
that
\[
\int_{\Omega}|\nabla\psi_{n}|^{2}=\left(  \int_{\Omega}a(x)\psi_{n}%
^{\;q+1}\right)  \Vert u_{n}\Vert^{q-1}+\alpha_{n}\int_{\partial\Omega}%
\psi_{n}^{\;2}\rightarrow0,
\]
so that $\psi_{n}\rightarrow\psi_{\infty}$ in $H^{1}(\Omega)$ and
$\psi_{\infty}$ is a positive constant.

On the other hand, since $u_{n}$ solves $(P_{\alpha_{n}})$ we have that
$\int_{\Omega}a(x)u_{n}^{\;q}=-\alpha_{n}\int_{\partial\Omega}u_{n}$. Dividing
by $\Vert u_{n}\Vert^{q}$ we obtain $\int_{\Omega}a(x)\psi_{n}^{\;q}%
=-\alpha_{n}\left(  \int_{\partial\Omega}\psi_{n}\right)  \Vert u_{n}%
\Vert^{1-q}$, and since
%\mpfns{modified}%
$\alpha_{n}\Vert u_{n}\Vert^{1-q}=\Vert w_{n}\Vert^{1-q}\rightarrow0$, we find
that $\int_{\Omega}a(x)\psi_{\infty}^{\;q}=0$. But $\psi_{\infty}$ is a
positive constant, so $\int_{\Omega}a=0$, contradicting $(A.0)$. \qed \newline

%\item Assume \marginpar{\scriptsize the proof is to be modified}by contradiction that we can choose $\alpha_{n}\rightarrow0^{+}$ and a nontrivial solution $w_{n}$ of $(R_{\alpha_{n}})$ such that $\Vert w_{n}\Vert_{C^{1}(\overline{\Omega})}\rightarrow0$ but \begin{equation} \Vert v_{n}\Vert_{C^{1}(\overline{\Omega})}\rightarrow0,\label{v1to0}% \end{equation}
%where $v_{n}=\alpha_{n}^{\;\frac{1}{q-1}}w_{n}$. Recall that $v_{n}$ satisfies \eqref{reducevn}. As in the proof of Proposition \ref{prop:bfzero}, conditions $(A.1)$ and $(A.2)$ allow us to assume that $v_{n}\not \equiv 0$ in some connected component $\Omega^{\prime}$ of $\Omega_{+}^{a}$. %Indeed, we observe that %\begin{align*}
%0<\int_{\Omega}|\nabla v_{n}|^{2}  &   =\int_{\Omega}a(x)v_{n}^{q+1}+\alpha
%_{n}\int_{\partial\Omega}b(x)v_{n}^{2}\\
%&  \leq\int_{\mathrm{supp}\,a^{+}}a^{+}(x)v_{n}^{q+1}+\alpha_{n}\int_{\Gamma
%_{+}^{b}}b(x)v_{n}^{2}\\
%&  =\int_{\Omega_{+}^{a}}a^{+}(x)v_{n}^{q+1}+\alpha_{n}\int_{\Gamma_{+}^{b}}%
%b^{+}(x)v_{n}^{2}.
%\end{align*}
%So, it should hold that $v_{n}\not \equiv 0$ in $\Omega_{+}^{a}$, in view of $(A.3)$. So, by using $(A.1)$, the desired assertion follows.
%Finally, \cite[Lemma 2.2]{KRQU16} is applicable to $v_{n}$ and $\Omega^{\prime}$. However, this contradicts \eqref{v1to0}. \qed

%\subsection{Local bifurcation analysis for $(R_{\alpha})$ via a
%Lyapunov-Schmidt type reduction}

We discuss now bifurcation of nontrivial solutions of $(R_{\alpha})$ from
$\Gamma_{1}$. To this end, we apply a Lyapunov-Schmidt type reduction. Let
\[
X_{2}:=\left\{  \psi\in L^{2}(\Omega):\int_{\Omega}\psi=0\right\}  .
\]
We decompose $w\in L^{2}(\Omega)$ as $w=t+\psi\in\mathbb{R}\oplus X_{2}$,
where $\psi:=Q[w]=w-t$ with $t:=(1/|\Omega|)\int_{\Omega}w$. By using the
projection $Q$ of $L^{2}(\Omega)$ into $X_{2}$, $(R_{\alpha})$ is reduced to
the following equations:
%we deduce that if
%$u$ is a nonnegative solution of $(R_{\alpha})$, then%
%\marginpar{{\tiny maybe we can avoid }$w${\tiny , and just use }$u=t+\psi
%${\tiny , as above}} \mpfns{{\bf agree. $u$ put back}}%
\begin{align*}
&  Q[-\Delta w-\alpha a(x)w^{q}]=0,\quad\partial_{\nu}w=\alpha w\quad
\mbox{on}\quad\partial\Omega,\\
&  (1-Q)[-\Delta w-\alpha a(x)w^{q}]=0.
\end{align*}
By direct calculations, with $w=t+\psi$, it follows that
%Note that $(P_\alpha)$ around the trivial solution $(\alpha, u) = (0, c_a)$ is reduced to the conditions%
\begin{align}
&
\begin{cases}
-\Delta\psi+\frac{\alpha}{|\Omega|}\int_{\partial\Omega}(t+\psi)=Q\left[
\alpha a(x)(t+\psi)^{q}\right]  & \mbox{in}\quad\Omega,\\
\partial_{\nu}\psi=\alpha(t+\psi) & \mbox{on}\quad\partial\Omega,
\end{cases}
\label{eq:Q}\\
&  \alpha\left(  \int_{\Omega}a(x)(t+\psi)^{q}+\int_{\partial\Omega}%
(t+\psi)\right)  =0. \label{eq:1-Q}%
\end{align}
First, we solve \eqref{eq:Q} around $(\alpha,t,\psi)=(0,c_{a},0)$. Let
%\marginpar{added}\mpfns{agree}
$r>N$ and%
\begin{align*}
&  W_{1}:=\left\{  \psi\in W^{2,r}(\Omega):\int_{\Omega}\psi=0\right\}  ,\\
&  Z_{1}:=\left\{  (g_{1},g_{2})\in L^{r}(\Omega)\times W^{1-\frac{1}{r}%
,r}(\partial\Omega):\int_{\Omega}g_{1}+\int_{\partial\Omega}g_{2}=0\right\}  .
\end{align*}
Let $B_{\delta}\subset W_{1}$ be a ball centered at the origin with radius
$\delta>0$. For a constant $c>0$, we define the nonlinear mapping
$\mathcal{F}:\mathbb{R}\times(c-\delta,c+\delta)\times B_{\delta}\rightarrow
Z_{1}$ by
\begin{align}
&  \mathcal{F}(\alpha,t,\psi)\nonumber\label{LS:F}\\
&  :=\left(  -\Delta\psi+\frac{\alpha}{|\Omega|}\int_{\partial\Omega}%
(t+\psi)-Q\left[  \alpha a(x)(t+\psi)^{q}\right]  , \partial_{\nu}\psi
-\alpha(t+\psi)\right)  .
\end{align}
Indeed, this is well defined, since $\int_{\Omega}Q\left[  \alpha
a(x)(t+\psi)^{q}\right]  =0$.
%For \marginpar{\textit{modified}}$c>0$,
Then, the Fr\'{e}chet derivative $\mathcal{F}_{\psi}(0,c,0):W_{1}\rightarrow
Z_{1}$ with respect to $\psi$ is given by $\mathcal{F}_{\psi}(0,c,0)\psi
=\left(  -\Delta\psi, \partial_{\nu}\psi\right)  $, and thus, it is a
homeomorphism. So, the implicit function theorem applies, and the equation
$\mathcal{F}(\alpha,t,\psi)=0$ is uniquely solvable around $(0,c,0)$ by some
$\psi=\psi(\alpha,t)$ satisfying $\psi(0,c)=0$.
%
%	because we can check by a simple calculation that \begin{align*}
%	\int_\Omega \left( -\Delta \psi + \frac{\alpha}{|\Omega|}\int_{\partial \Omega}b(x)(t+\psi) - Q\left[ \alpha a(x) (t+\psi)^q \right] \right) + \int_{\partial \Omega} \left( \frac{\partial \psi}{\partial \nu} - \alpha b(x) (t+\psi) \right)
%	= 0.
%	\end{align*}
%	Indeed,  is uniquely solvable with respect to $\psi$ at $(0, c_a, 0)$ in the space
%	$\mathbb{R} \times \mathbb{R} \times W_1 := \{ \psi \in W^{2,r}(\Omega) : \int_\Omega \psi = 0 \}$. The unique solution is then given by
%	\begin{align*}
%	& \psi = \psi (\alpha, t) \quad\mbox{for} \quad (\alpha, t) \simeq (0, c_a), \\
%	& \psi (0, c_a)=0.
%	\end{align*}
Plugging $\psi(\alpha,t)$ into \eqref{eq:1-Q}, we obtain the bifurcation
equation
\begin{equation}
\alpha\left(  \int_{\Omega}a(x)(t+\psi(\alpha,t))^{q}+\int_{\partial\Omega
}(t+\psi(\alpha,t))\right)  =0. \label{eq:1-Q:t}%
\end{equation}

Summing up, solving $(R_{\alpha})$ around $(\alpha, w)=(0,c)$ reduces to the
solvability of the equation
\begin{align}
\label{LS:beq}\mathcal{G}(\alpha, t) := \int_{\Omega}a(x) (t+ \psi(\alpha,
t))^{q} + \int_{\partial\Omega} (t + \psi(\alpha, t)) = 0,
\end{align}
around $(\alpha, t)=(0, c)$ (note that $\alpha=0$ in \eqref{eq:1-Q:t} yields
the trivial solution $(\alpha, w) = (0,c)$).
%
%Indeed, we see that if a positive function $t=t(\alpha)$ solves $\mathcal{G}(\alpha, t)=0$, then putting
%\begin{align}
%\label{umu}u(\alpha) := t(\alpha) + \psi(\alpha, t(\alpha)),
%\end{align}
%we infer from \eqref{eq:Q} that
%\begin{align*}%
%\begin{cases}
%-\Delta u(\alpha) + \frac{\alpha}{|\Omega|}\int_{\partial\Omega}b(x) u(\alpha) = Q
%\left[  \alpha a(x) u(\alpha)^{q} \right]  & \mbox{in} \quad\Omega,\\
%\frac{\partial u(\alpha)}{\partial\nu} = \alpha b(x) u(\alpha) & \mbox{on}
%\quad\partial\Omega.
%\end{cases}
%\end{align*}
%Moreover, we observe from the assertion $\mathcal{G}(\alpha, t(\alpha))=0$ that
%\begin{align*}
%Q \left[  \alpha a(x) u(\alpha)^{q} \right]   &  = \alpha a(x) u(\alpha)^{q} - %\frac{\alpha
%}{|\Omega|} \int_{\Omega}a(x) u(\alpha)^{q}\\
%&  = \alpha a(x) u(\alpha)^{q} + \frac{\alpha}{|\Omega|} \int_{\partial\Omega}
%b(x)u(\alpha).
%\end{align*}
%Hence, we obtain that
%\begin{align*}%
%\begin{cases}
%-\Delta u(\alpha) = \alpha a(x) u(\alpha)^{q} & \mbox{in} \quad\Omega,\\
%\frac{\partial u(\alpha)}{\partial\nu} = \alpha b(x)u(\alpha) & \mbox{on} \quad
%\partial\Omega,
%\end{cases}
%\end{align*}
%as desired.
%We remark that $\alpha = 0$ is a solution of \eqref{eq:1-Q:t}, and in this case, plugging $\alpha = 0$ into \eqref{eq:Q}, we find that $\psi = 0$, so $u$ is a constant. This means that $(\alpha, u)$ is on $\Gamma_1$.

In the sequel we prove that under $(A.0)$ a certain $C^{\infty}$ mapping
$\alpha\mapsto t(\alpha)$ uniquely solves \eqref{LS:beq} around $(0,c_{a})$.
Conversely, we show that, besides $(0,0)$, this is the only bifurcation point
in $\Gamma_{1}$ for solutions of $(R_{\alpha})$. More generally, we prove that
$(0,0)$ and $(0,c_{a})$ are the only possible limits for a sequence
$(\alpha_{n},w_{n})$ with $\alpha_{n}\rightarrow0^{+}$ and $w_{n}$ solving
$(R_{\alpha_{n}})$.

%------------------------------------

%-----------------------------------

\begin{figure}[tbh]
\begin{center}
\includegraphics[scale=0.2]{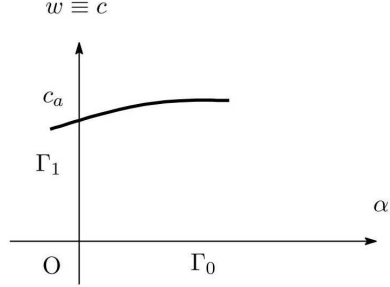}
\end{center}
\caption{Bifurcating positive solutions of $(R_{\alpha})$ at $(0, c_{a})$.}%
\label{fig18_0112a}%
\end{figure}

\begin{prop}
[Bifurcation from $\Gamma_{1}$]\label{prop:bG00} Assume $(A.0)$. Then:

\begin{enumerate}
\item $(R_{\alpha})$ has solutions $w=w(\alpha)\in\mathcal{P}^{\circ}$
bifurcating from $\Gamma_{1}$ at $(0,c_{a})$, and such that
%\mpfns{{\bf modified}}
$\alpha\mapsto w(\alpha)=t(\alpha)+\psi(\alpha,t(\alpha))$ is $C^{\infty}$
from $(-\alpha_{0},\alpha_{0})$ into $W^{2,r}(\Omega)$ for some $\alpha_{0}%
>0$, and $t(0)=c_{a}$, where $w=t+\psi$ is the decomposition as above.
Moreover, if $(\alpha,w)$ is a solution of $(R_{\alpha})$ around $(0,c_{a})$
in $\mathbb{R}\times C^{1}(\overline{\Omega})$, then $w=w(\alpha)$ for some
$\alpha$, see Figure \ref{fig18_0112a}.

\item Let $\alpha_{n}\rightarrow0^{+}$ and $w_{n}$ be nontrivial solutions of
$(R_{\alpha_{n}})$. Then, up to a subsequence, we have either $w_{n}%
\rightarrow0$ or $w_{n}\rightarrow c_{a}$ in $C(\overline{\Omega})$.
\end{enumerate}
\end{prop}

%-------------------------------------------------------------

%\mpfns{a remark added}

%Figure 5

%	% \begin{figure}[H]

\noindent\textit{Proof.}

\begin{enumerate}
\item First of all, let us observe that once we get positive solutions
bifurcating from $\Gamma_{1}$ at $(0, c_{a})$ in $\mathbb{R}\times
W^{2,r}(\Omega)$, these ones are in $\mathcal{P}^{\circ}$, since
$W^{2,r}(\Omega)\subset C^{1}(\overline{\Omega})$ and $c_{a}$ is a positive constant.

We use the implicit function theorem to analyze the reduced bifurcation
equation $\mathcal{G}(\alpha, t)=0$ around $(0, c_{a})$. Note from
\eqref{def:ca} that $\mathcal{G}(0, c_{a}) = 0$. Differentiating $\mathcal{G}$
with respect to $t$ yields
\begin{align*}
\frac{\partial\mathcal{G}}{\partial t}(\alpha, t) = \int_{\Omega}%
a(x)q(t+\psi)^{q-1} \left(  1 + \frac{\partial\psi}{\partial t} \right)  +
\int_{\partial\Omega} \left(  1 + \frac{\partial\psi}{\partial t} \right)  .
\end{align*}
From \eqref{eq:Q}, we see that $\frac{\partial\psi}{\partial t}(0, c_{a})=0$,
so that $\frac{\partial\mathcal{G}}{\partial t}(0, c_{a}) = \int_{\Omega}a(x)q
c_{a}^{\; q-1} + |\partial\Omega|$. Using \eqref{def:ca}, it follows that
$\frac{\partial\mathcal{G}}{\partial t}(0, c_{a}) = (1-q)|\partial\Omega| >
0$. The implicit function theorem is now applicable, and then, we obtain that
for $(\alpha, t)\simeq(0, c_{a})$,
\begin{align}
\mathcal{G}(\alpha, t)=0 \ \Longleftrightarrow\ t=t(\alpha) \ \mbox{with}
\ \ t(0)=c_{a}. \label{implicit}%
\end{align}
%Thus, $u(\alpha)$ defined by \eqref{umu} for $\alpha\simeq0$ are the desired positive solutions,
%and there are no other nontrivial solutions of $(R_{\alpha})$ around $(0, c_{a})$.
Finally, the assertion that $\alpha\mapsto w(\alpha)=t(\alpha)+\psi
(\alpha,t(\alpha))$ is $C^{\infty}$ follows from the well known regularity
argument for the implicit function theorem.

The uniqueness assertion can be verified in a similar way as in the proof of
Proposition \ref{last}(ii).

\item Since $w_{n}$ solves $(R_{\alpha_{n}})$ with $\alpha_{n}>0$, we know by
Proposition \ref{cor:bounds} that $\{w_{n}\}$ is bounded in $C^{1}%
(\overline{\Omega})$, and consequently in $H^{1}(\Omega)$. Thus, up to a
subsequence, we have $w_{n}\rightharpoonup w$ in $H^{1}(\Omega)$ and
$w_{n}\rightarrow w$ in $L^{t}(\Omega)$ for $1\leq t<2^{\ast}$, and in
$L^{2}(\partial\Omega)$. Taking the limit as $n\rightarrow\infty$ in the weak
formulation of $(R_{\alpha_{n}})$ we see that $w_{n}\rightarrow w$ in
$C(\overline{\Omega})$ and $w$ is a nonnegative constant. Moreover, from
$\int_{\Omega}a(x)w_{n}^{q}+\int_{\partial\Omega}w_{n}=0$ we obtain
$w^{q}\left(  \int_{\Omega}a+w^{1-q}|\partial\Omega|\right)  =0$, so either
$w=0$ or $w=c_{a}$. \qed

\end{enumerate}

\begin{rem}
\label{rpq} Proposition \ref{prop:bG00}(i) can be formulated in a more general
setting as follows:
%Instead of \eqref{LS:F}, we reset the mapping $\mathcal{F}:\mathbb{R}\times(-\delta, 1)\times(c-\delta,c+\delta)\times B_{\delta} \rightarrow Z_{1}$ by
%\begin{align*}
%&  \mathcal{F}(\alpha,q,t,\psi)\\
%&  =\left(  -\Delta\psi+\frac{\alpha}{|\Omega|} \int_{\partial\Omega}(t+\psi)-Q\left[  \alpha a(x) (t+\psi)^{q}\right]  , \frac{\partial\psi}{\partial\nu}-\alpha( t+\psi)\right)  .
%\end{align*}
Assume $(A.0)$, $q_{0}\in\lbrack0,1)$, and let $c_{a,q_{0}}$ be $c_{a}$ with
$q=q_{0}$. Then
%we reduce $(R_{\alpha})$ around $(\alpha,q,w)=(0,q_{0},c)$ to \eqref{LS:beq} replaced by
%\[
%\mathcal{G}(\alpha,q,t)=\int_{\Omega}a(x) (t+ \psi(\alpha,q,t))^{q} + \int_{\partial\Omega} (t + \psi(\alpha,q,t)) = 0
%\]
%around $(\alpha,q,t)=(0,q_{0},c)$. Consequently, under the assumptions of Proposition \ref{prop:bG00},
$(R_{\alpha})$ has, around $(\alpha,q,w)=(0,q_{0},c_{a,q_{0}})$, exactly one
solution $w=w(\alpha,q)\in\mathcal{P}^{\circ}$ parametrized by $(\alpha,q)$,
and such that $(\alpha,q)\mapsto w(\alpha,q)=t(\alpha,q)+\psi(\alpha
,q,t(\alpha,q))$ is $C^{\infty}$ from $(-\alpha_{0},\alpha_{0})\times
(q_{0}-\delta_{0},q_{0}+\delta_{0})$ into $W^{2,r}(\Omega)$ for some
$\alpha_{0},\delta_{0}>0$, and $t(0,q_{0})=c_{a,q_{0}}$.
\end{rem}

%	% \begin{figure}[H]
%\begin{figure}[tbh]
%\begin{center}
%\includegraphics[scale=0.18]{figLS2.eps}
%\end{center}
%\caption{The solutions $w=t+\psi(\alpha,q,t)\in\mathcal{P}^{\circ}$ of $(R_{\alpha})$ parametrized by $(\alpha,q)$.}%
%\label{figLS2}%
%\end{figure}

As a corollary of Theorem \ref{thm:curve} and Propositions \ref{cor:bounds},
\ref{prop:C1C2} and \ref{prop:bG00}, we obtain the following exact
multiplicity result for $(R_{\alpha})$:

\begin{cor}
[Exact multiplicity for $(R_{\alpha})$]\label{cor:uqla} Assume $(A.0)$ and
$0<\delta<c_{a}$. Then there exists $\alpha_{\delta}>0$ such that, for each
$0<\alpha<\alpha_{\delta}$:

\begin{enumerate}
\item $(R_{\alpha})$ has a unique solution $w$ satisfying $\Vert w\Vert
_{C^{1}(\overline{\Omega})}>\delta$. Moreover, $w=w(\alpha)\in\mathcal{P}%
^{\circ}$ from Proposition \ref{prop:bG00}(i).

\item If we assume, in addition, $(A.1)$ and $q \in\mathcal{A}_{N}$, then
$(R_{\alpha})$ has a unique nontrivial solution $w$ satisfying $\Vert
w\Vert_{C^{1}(\overline{\Omega})}<\delta$, namely, $w=\alpha^{\frac{1}{1-q}%
}u_{\alpha}\in\mathcal{P}^{\circ}$, where $u_{\alpha}$ is given by Theorem
\ref{thm:curve}.
\end{enumerate}
\end{cor}

\noindent\textit{Proof.} The first item follows promptly from Propositions
\ref{cor:bounds} and \ref{prop:bG00}. We prove now the second item. By Theorem
\ref{thm:curve} we know that $\alpha^{\frac{1}{1-q}}u_{\alpha}$ solves
$(R_{\alpha})$. We claim that it is the only solution of $(R_{\alpha})$
converging to $0$ in $C^{1}(\overline{\Omega})$ as $\alpha\rightarrow0^{+}$.
Indeed, by Proposition \ref{prop:C1C2}, if $w_{n}$ is such a solution then
$\{u_{n}:=\alpha_{n}^{-\frac{1}{1-q}}w_{n}\}$ remains bounded in
$C^{1}(\overline{\Omega})$. Hence, $u_{n}\rightarrow u_{N}$ in $C^{1}%
(\overline{\Omega})$ by elliptic regularity, Lemma \ref{prop:glb}, and the
condition $q\in\mathcal{A}_{N}$. By Theorem \ref{thm:curve}(ii) we infer that
for $n$ large enough $u_{n}=u_{\alpha}$ for some $\alpha>0$. The proof is now
complete. \qed\newline

We end this section with the corresponding exact multiplicity result for
$(P_{\alpha})$, which follows from Corollary \ref{cor:uqla}:
%----------------------------------------------------------------
%Proposition

\begin{theorem}
[Exact multiplicity for $(P_{\alpha})$]\label{prop:exactm} Assume $(A.0)$,
$(A.1)$, and $q\in\mathcal{A}_{N}$. Then there exists $\alpha_{0}>0$ such that
$(P_{\alpha})$ has exactly two nontrivial solutions $u_{1,\alpha},u_{2,\alpha
}$ for $0<\alpha<\alpha_{0}$. Moreover, $u_{1,\alpha},u_{2,\alpha}%
\in\mathcal{P}^{\circ}$ and $u_{1,\alpha}<u_{2,\alpha}$ on $\overline{\Omega}$.
\end{theorem}

\section{A topological bifurcation approach to $(R_{\alpha})$}

\label{sec:topoba}

The proof of Theorem \ref{maint} is based on a bifurcation approach via a
\textit{regularization} scheme, which analyzes the structure of the solutions
set of $(R_{\alpha})$. More precisely, we study how the bifurcation curve
obtained by Proposition \ref{prop:bG00} behaves globally in $\alpha>0$.
%\marginpar{modified}

Introducing a new parameter $\varepsilon\in(0,1]$, we consider
\[%
\begin{cases}
-\Delta w=\alpha a(x)(w+\varepsilon)^{q-1}w\ =\alpha a(x)\left(  \frac
{w}{w+\varepsilon}\right)  ^{1-q}w^{q} & \mbox{ in }\Omega,\\
w \geq0 & \mbox{ in }\Omega,\\
\partial_{\nu} w=\alpha w & \mbox{ on }\partial\Omega.
\end{cases}
\leqno{(R_{\alpha}^\varepsilon)}
\]
Note that any nontrivial solution of $(R_{\alpha}^{\varepsilon})$ belongs to
$\mathcal{P}^{\circ}$, since $s\mapsto(s+\varepsilon)^{q-1}s$ is $C^{1}$ in
$[0,\infty)$, and consequently the strong maximum principle and Hopf's lemma apply.

%---------------------------------------------

%\marginpar{{\tiny subsection removed, and phrase added}}
We start with some preliminary results, namely, the counterparts of
Propositions \ref{cor:bounds}, \ref{alpha2} and \ref{prop:bG00}(ii) for
$(R_{\alpha}^{\varepsilon})$.
%The following three results are the counterparts of Propositions \ref{prop:bpc}
%{cor:bounds} %and \ref{prop:uppblam}(i) and Corollary \ref{cor:bfi}, respectively, for $(R_{\alpha}^{\varepsilon})$.

We establish an \textit{a priori} estimate in $C^{1}(\overline{\Omega})$ for
solutions in $\mathcal{P}^{\circ}$ of $(R_{\alpha}^{\varepsilon})$, i.e. the
counterpart of Proposition \ref{cor:bounds}:

\begin{prop}
\label{cor:bfie}
%Assume $\left(  A.3\right)  $.
Let $\overline{\alpha}, \overline{\varepsilon}>0$. Then there exists
$\overline{C}>0$ such that $\Vert w\Vert_{C^{1}(\overline{\Omega})}%
\leq\overline{C}$ for every solution $w\in\mathcal{P}^{\circ}$ of $(R_{\alpha
}^{\varepsilon})$ with $\alpha\in(0,\overline{\alpha}]$ and $\varepsilon
\in(0,\overline{\varepsilon}]$.
\end{prop}

\noindent\textit{Proof.} Assume by contradiction that $w_{n}\in\mathcal{P}%
^{\circ}$ is a solution of $(R_{\alpha_{n}}^{\varepsilon_{n}})$ such that
$0<\alpha_{n}\rightarrow\alpha_{\infty}\in\lbrack0,\overline{\alpha}]$,
$\varepsilon_{n}\in(0,\overline{\varepsilon}]$ but $\Vert w_{n}\Vert
_{C^{1}(\overline{\Omega})}\rightarrow\infty$. We can then argue as in the
proof of Propositions \ref{cor:bounds}, with $w_{n}^{q}$ replaced by $\left(
\frac{w_{n}}{w_{n}+\varepsilon_{n}}\right)  ^{1-q}w_{n}^{q}$, and notice that
$\left\vert \frac{w_{n}}{w_{n}+\varepsilon_{n}}\right\vert \leq1$ for
%\marginpar{{\tiny subsection removed}}
$n\geq1$. \qed\newline

The next proposition is the counterpart of Proposition \ref{prop:bG00}(ii).

\begin{prop}
\label{prop:bpc:e} Assume $(A.0)$ and $\varepsilon\in(0,c_{a})$. If $w_{n}%
\in\mathcal{P}^{\circ}$ are solutions of $(R_{\alpha_{n}}^{\varepsilon})$ with
$\alpha_{n}\rightarrow0^{+}$ then, up to a subsequence, we have either
$w_{n}\rightarrow0$ or $w_{n}\rightarrow c_{a}-\varepsilon$ in $C(\overline
{\Omega})$.
\end{prop}

\noindent\textit{Proof.} We use Proposition \ref{cor:bfie} and argue as in the
proof of Proposition \ref{prop:bG00}(ii) to deduce that, up to a subsequence,
$w_{n}\rightarrow c$ in $C(\overline{\Omega})$, where $c$ is a nonnegative
constant such that $0=\left(  c+\varepsilon\right)  ^{q-1}c\int_{\Omega
}a+c|\partial\Omega|$. The desired conclusion thus follows. \qed\newline
%Since $w_{n}$ satisfies $(P_{\alpha_{n}, \varepsilon})$,it follows that
%	\[
%	\alpha_n \left( \int_\Omega a (u_n+\varepsilon)^{q-1}u_n + \int_{\partial \Omega} b u_n \right)
%	= 0,
%	\]
%	so we have%
%\begin{align}
%\label{bfc:une}\int_{\Omega}a (w_{n}+\varepsilon)^{q-1}w_{n} + \int_{\partial\Omega} b w_{n} = 0.
%\end{align}
%Passing to the limit as $n\to\infty$, we have
%\begin{align*}
%(c+\varepsilon)^{q-1}c \int_{\Omega}a + c \int_{\partial\Omega} b = 0.
%\end{align*}
%The desired conclusion follows. \qed \newline
%	The following result is the counterpart of Proposition \ref{prop:bfi0} for $(P_{\lambda, \varepsilon})$.
%	\begin{prop} \label{prop:bfi0:e}
%	Assume $\int_{\partial \Omega} b \neq 0$. Let $\varepsilon \in (0,1]$. Then, there is no bifurcation of solutions in $\mathcal{P}^\circ$ from infinity for $(P_{\lambda, \varepsilon})$ at $\lambda = 0$.
%	\end{prop}
%	\begin{proof}
%	We only note that $(u+\varepsilon)^{q-1}u = \left( \frac{u}{u+\varepsilon} \right)^{1-q} u^q$. Thus, it is straightforward.
%	\end{proof}

%	The following result is the counterpart of Proposition \ref{bfin0} for $(P_{\lambda, \varepsilon})$.

%	\begin{prop} \label{bfin0:e}
%	Assume $\int_{\partial \Omega} b > 0$. Let $\varepsilon \in (0,1]$. Then, there is no bifurcation of solutions in $\mathcal{P}^\circ$ from infinity for $(P_{\lambda, \varepsilon})$ at any $\lambda > 0$.
%	\end{prop}
%	\begin{proof}
%	It is straightforward.
%	\end{proof}

%\marginpar{removed}

%---------------- bounds of norm

%----------------------- bounds of \alpha

Next we establish an \textit{a prori} upper bound of $\alpha>0$ for the
existence of a solution in $\mathcal{P}^{\circ}$ of $(R_{\alpha}^{\varepsilon
})$. Using \eqref{changeofv}, we reduce $(R_{\alpha}^{\varepsilon})$ to the
problem
\[%
\begin{cases}
-\Delta u=a(x)\left(  \frac{\alpha^{\frac{1}{1-q}}u}{\alpha^{\frac{1}{1-q}%
}u+\varepsilon}\right)  ^{1-q}u^{q} & \mbox{ in }\Omega,\\
u \geq0 & \mbox{ in }\Omega,\\
\partial_{\nu} u=\alpha u & \mbox{ on }\partial\Omega.
\end{cases}
\leqno{(P_{\alpha}^\varepsilon)}
\]
We remark that, as long as $\alpha>0$, $w$ solves $(R_{\alpha}^{\varepsilon})$
if and only if $u=\alpha^{-\frac{1}{1-q}}w$ solves $(P_{\alpha}^{\varepsilon
})$. So, it suffices to establish the upper bound for $(P_{\alpha
}^{\varepsilon})$.

The following result is the counterpart of Proposition \ref{alpha2}.

\begin{prop}
\label{prop:nbddalph} Assume $(A.3)$. Then, for any $q\in\lbrack0,1)$ there
exist $\underline{\alpha},\overline{\varepsilon}>0$ such that $(R_{\alpha
}^{\varepsilon})$ or $(P_{\alpha}^{\varepsilon})$ has no solutions in
$\mathcal{P}^{\circ}$ for any $\alpha>\underline{\alpha}$ and $0<\varepsilon
\leq\overline{\varepsilon}$.
\end{prop}

\noindent\textit{Proof.} It suffices to consider the case $(P_{\alpha
}^{\varepsilon})$, taking \eqref{changeofv} into account. Let $u\in
\mathcal{P}^{\circ}$ be a solution of $(P_{\alpha}^{\varepsilon})$ with
$\alpha>0$ and $\varepsilon>0$. Then, Green's formula yields
\[
\int_{\Omega}a(x)\left(  \frac{\alpha^{\frac{1}{1-q}}u}{\alpha^{\frac{1}{1-q}%
}u+\varepsilon}\right)  ^{1-q}=\int_{\Omega}\frac{-\Delta u}{u^{q}}%
<-\alpha\int_{\partial\Omega}u^{1-q}.
\]
It follows that
%\marginpar{\scriptsize replaced by $\int_\Omega a^-$}%
\begin{equation}
\alpha\int_{\partial\Omega}u^{1-q}<\int_{\Omega}a^{-}. \label{boundCa}%
\end{equation}

The rest of the proof proceeds in a similar manner as in the proof of
Proposition \ref{alpha2}. Indeed, let $\beta,\varepsilon_{0}>0$. In place of
$(Q_{q,a})$, we consider the following concave mixed problem:
\begin{equation}%
\begin{cases}
-\Delta u=a(x)f_{\beta,\varepsilon_{0}}(u) & \mbox{in $D$},\\
\partial_{\nu}u=0 & \mbox{on $\Gamma_{\partial \Omega}$},\\
u=0 & \mbox{on $\Gamma_{\Omega}$},
\end{cases}
\label{cp}%
\end{equation}
where
\[
f_{\beta,\varepsilon}(s):=\left(  \frac{\beta^{\frac{1}{1-q}}s}{\beta
^{\frac{1}{1-q}}s+\varepsilon}\right)  ^{1-q}s^{q},\quad s\geq0.
\]
Note that $s\mapsto f_{\beta,\varepsilon}(s)/s$ is decreasing for $s>0$. Since
$f_{\beta,\varepsilon}(s)$ is increasing with respect to $\beta>0$ and
decreasing with respect to $\varepsilon>0$ for every $s>0$, $u$ is a
supersolution of \eqref{cp} for $\alpha\geq\beta$ and $0<\varepsilon
\leq\varepsilon_{0}$. Consequently, given $\beta> 0$, we can choose
$\varepsilon_{0} > 0$ small enough such that, denoting by $u_{\beta
,\varepsilon_{0}}$ the unique (by Lemma \ref{lem:ump}) positive solution of
\eqref{cp} satisfying $u_{\beta,\varepsilon_{0}}>0$ in $D\cup\Gamma
_{\partial\Omega}$ (which exists, as in Proposition \ref{alpha2}), we have
that
\begin{equation}
u\geq u_{\beta,\varepsilon_{0}}\quad\mbox{on $\overline{D}$} \label{u:lbD}%
\end{equation}
for $\alpha\geq\beta$ and $0<\varepsilon\leq\varepsilon_{0}$. Combining
\eqref{boundCa} with \eqref{u:lbD} provides the desired conclusion.
\qed\newline

Next, under $(A.0)$, we will prove the existence of positive solutions of
$(R_{\alpha}^{\varepsilon})$ bifurcating from $\Gamma_{0}$. To obtain
bifurcation points from $\Gamma_{0}$ for positive solutions, we consider the
linearized eigenvalue problem at $w=0$:
\begin{equation}%
\begin{cases}
-\Delta\phi=\alpha a(x)\varepsilon^{q-1}\phi & \mbox{ in }\Omega,\\
\partial_{\nu}\phi=\alpha\phi & \mbox{ on }\partial\Omega.
\end{cases}
\label{lepro}%
\end{equation}
Since $(A.0)$ implies that $\int_{\Omega}a\varepsilon^{q-1}+|\partial
\Omega|<0$ and $(a\varepsilon^{q-1})^{+}\not \equiv 0$ if $\varepsilon$ is
small enough, \eqref{lepro} has exactly two principal eigenvalues, namely,
$\alpha=0,\alpha_{1,\varepsilon}$, where $\alpha_{1,\varepsilon}>0$ and both
are simple. Moreover, \eqref{lepro} possesses positive eigenfunctions
$\phi_{0},\phi_{1,\varepsilon}$ associated to $0,\alpha_{1,\varepsilon}$,
respectively, where $\phi_{0}$ is a positive constant (see \cite[Theorem
2.1]{U10}).

Applying to both $(0,0)$ and $(\alpha_{1,\varepsilon}, 0)$ the local and
unilateral global bifurcation theory from simple eigenvalues
\cite{CR71,Ra71,LGbook}, we obtain two \textit{components} (\textit{i.e.},
nonempty, \textit{maximal} closed and connected subsets) $\gamma_{0,
\varepsilon}$, $\gamma_{1, \varepsilon}$ in $\mathbb{R} \times C^{1}%
(\overline{\Omega})$ of solutions of $(R_{\alpha}^{\varepsilon})$, containing
$(0,0)$ and $(\alpha_{1,\varepsilon}, 0)$, respectively. In addition,
$\gamma_{0, \varepsilon}, \gamma_{1,\varepsilon}$
%\mpfns{`maximal' deleted}
%are maximal (for inclusion) and
consist of solutions in $\mathcal{P}^{\circ}$ except $(0,0)$, $(\alpha
_{1,\varepsilon}, 0)$. Moreover, the set of nontrivial solutions of
$(R_{\alpha}^{\varepsilon})$ near $(0,0)$, $(\alpha_{1, \varepsilon}, 0)$ is
given exactly by $\gamma_{0, \varepsilon}$, $\gamma_{1, \varepsilon}$,
respectively, so $\Gamma_{1} \subset\gamma_{0, \varepsilon}$.

Based on the existence of $\gamma_{0,\varepsilon}$, $\gamma_{1,\varepsilon}$,
we state the main result of this section, which extends the local existence
and multiplicity result proved in \cite[Theorem 5.2, Proposition 7.4, Lemma
7.5]{CT14} by showing that $(R_{\alpha})$ has a subcontinuum of nontrivial
solutions for $\alpha>0$.

\begin{theorem}
\label{maint} Assume $(A.0)$, $(A.1)$, $(A.2)$, $(A.3)$ and $q\in
\mathcal{A}_{N}$. Then $(R_{\alpha})$ possesses a subcontinuum $\gamma_{0}$ in
$[0,\infty)\times C^{1}(\overline{\Omega})$ of solutions satisfying
\eqref{gcapGab} (see Figure \ref{fig18_0105b}). Moreover, the following
\textrm{three} assertions hold:

\begin{enumerate}
%\item $(\alpha, u) \in\gamma_{0} \setminus\{ (0,0), (0, c_a)\}$ implies that $\alpha> 0$ and $u\not \equiv 0$. In particular, $\gamma_{0}$ does not meet any point on $\Gamma_{0} \cup\Gamma_1$ except $(0,0), (0, c_a)$.

\item $u\in\mathcal{P}^{\circ}$ for $(\alpha,u)\in\gamma_{0}\setminus\{
(0,0)\}$.

\item There exists $\overline{\alpha}>0$ such that $(R_{\alpha})$ has exactly
two nontrivial solutions $w_{1,\alpha},w_{2,\alpha}$ for $0<\alpha<
\overline{\alpha}$, which satisfy $(\alpha,w_{1,\alpha})$, $(\alpha
,w_{2,\alpha})\in\gamma_{0}$, and $w_{1,\alpha}<w_{2,\alpha}$ on
$\overline{\Omega}$. Furthermore,
%In particular, $w_{2,\alpha}\in\mathcal{P}^{\circ}$.

\begin{enumerate}
\item $w_{1,\alpha}=\alpha^{\frac{1}{1-q}}u_{\alpha}$, where $u_{\alpha}$ is
given by Theorem \ref{thm:curve}(ii).

\item $w_{2,\alpha}=w(\alpha)$, where $w(\alpha)$ is given by Proposition
\ref{prop:bG00}.
%
%$\alpha \mapsto w_{2,\alpha}$ is a $C^{\infty}$ mapping from $[0, \overline{\alpha})$ to  $W^{2,r}(\Omega)$ such that $w_{2,0}=c_{a}$.
%In particular, $w(\alpha)\in\mathcal{P}^{\circ}$.
%Moreover, any nontrivial solution $w$ of $(R_{\alpha})$ around $(0,c_{a})$ in $\{ \alpha\geq0 \} \times C^{1}(\overline{\Omega})$ lies on $\gamma_{0}$.

\end{enumerate}

%\item

%\item Let $\overline{\alpha}>0$. Then, there exists $C>0$ such that $\|w\|_{C^{1}(\overline{\Omega})}\leq C$ for $(\alpha,w)\in\gamma_{0}$ with $\alpha\in[0,\overline{\alpha}]$.
%	\item For any $(\alpha, u)\in\gamma_{0}$ in a neighborhood of $(0, c_a)$, $u$
%	is a positive solution of $(R_{\alpha})$ in $\mathcal{P}^{\circ}$. Moreover, any
%	nontrivial solution $u$ of $(R_{\alpha})$, $\alpha\geq0$, around
%	$(0,c_a)$ lies on $\gamma_{0}$.

%	\item Assume $q\in\mathcal{A}$. Then, there exists $\overline{\alpha}>0$ such
%	that any nontrivial solution of $(R_{\alpha})$ for $0<\alpha
%	<\overline{\alpha}$ is in $\mathcal{P}^{\circ}$. Additionally

%	\marginpar{we do
%	not need Holder`s assumptions, right?\\ {\it\small right.
%	Thm\ref{maint} modified. }}

%	assume $a\in C^{\theta}(\overline
%	{\Omega})$ and $b\in C^{1+\theta}(\partial\Omega)$. Then, $(R_{\alpha})$
%	possesses \textbf{exactly two} nontrivial solutions $u_{1,\alpha
%	},u_{2,\alpha}$ for $0<\alpha<\overline{\alpha}$ (made smaller if necessary),
%	satisfying $0<u_{1,\alpha}<u_{2,\alpha}$ in $\overline{\Omega}$.

\item Let $\hat{\gamma}_{0}$ be the \textrm{component} of solutions of
$(R_{\alpha})$ in $[0,\infty)\times C^{1}(\overline{\Omega})$ that contains
$\gamma_{0}$. Then $\hat{\gamma}_{0}\setminus(\Gamma_{0}\cup\Gamma_{1})$ is
bounded in $[0,\infty)\times C^{1}(\overline{\Omega})$ (and in $[0,\infty
)\times W^{2,r}(\Omega)$, by elliptic regularity) and composed by solutions in
$\mathcal{P}^{\circ}$. In addition,
\begin{align}  \label{hatgam0}
(\hat{\gamma}_{0}\setminus(\Gamma_{0}\cup\Gamma_{1}))\cap\left\{
(\alpha,0):\alpha>0\right\}  =\emptyset.
\end{align}

%	\item Assume $(A.2^{\prime})$ and the condition $q\in\mathcal{A}$, and
%	additionally assume $b\geq0$. Then, $\gamma_{0}\setminus\{(0,0),(0,c_a)\}$
%	consists of positive solutions in $\mathcal{P}^{\circ}$.

\end{enumerate}
\end{theorem}

\noindent\textit{Proof.} First of all, by Proposition \ref{cor:Aal:ordered},
every nontrivial solution of $(R_{\alpha})$ lies in $\mathcal{P}^{\circ}$.

\begin{figure}[tbh]
\centerline{
\includegraphics[scale=0.18]{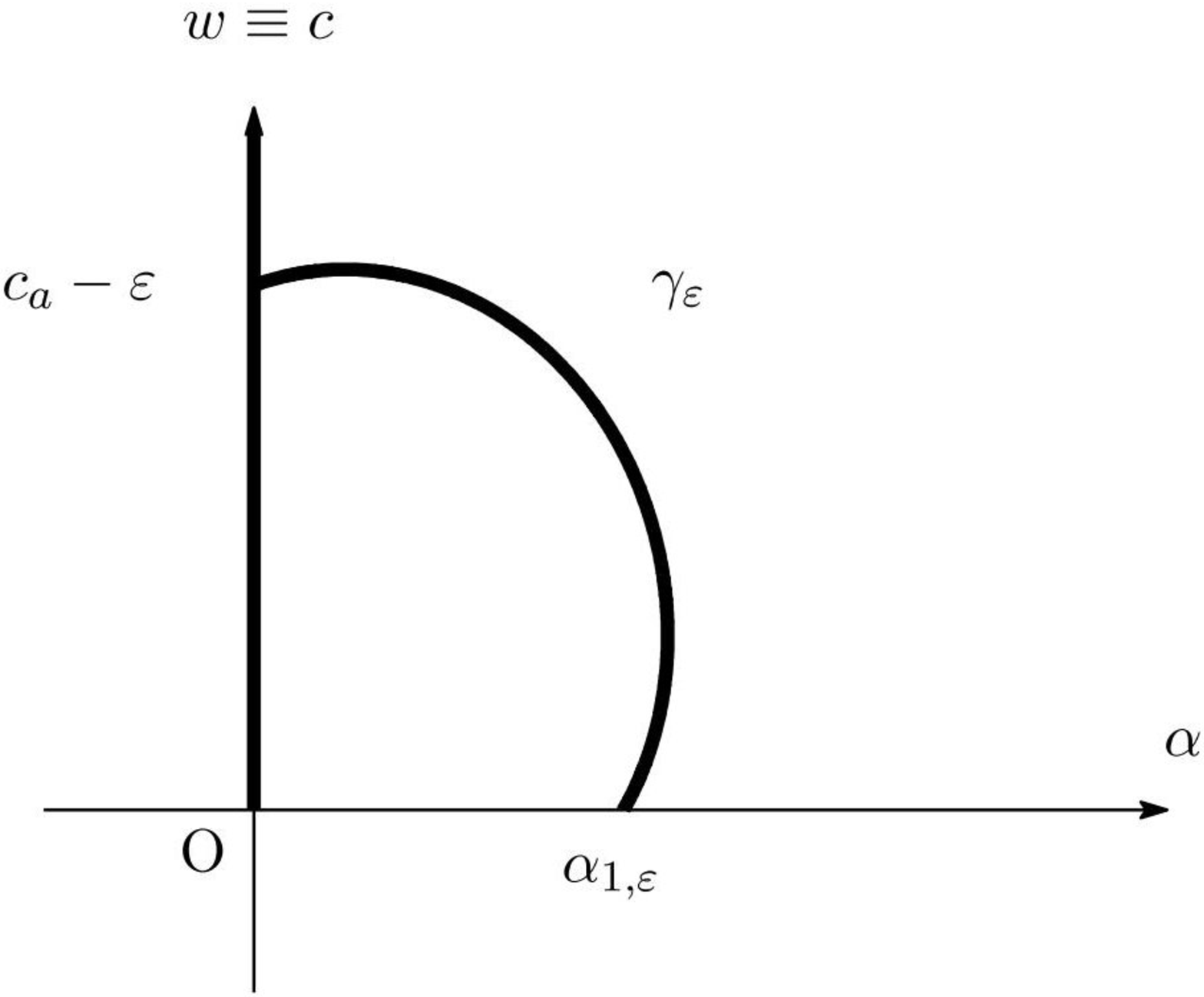} \hskip0.35cm
\includegraphics[scale=0.18] {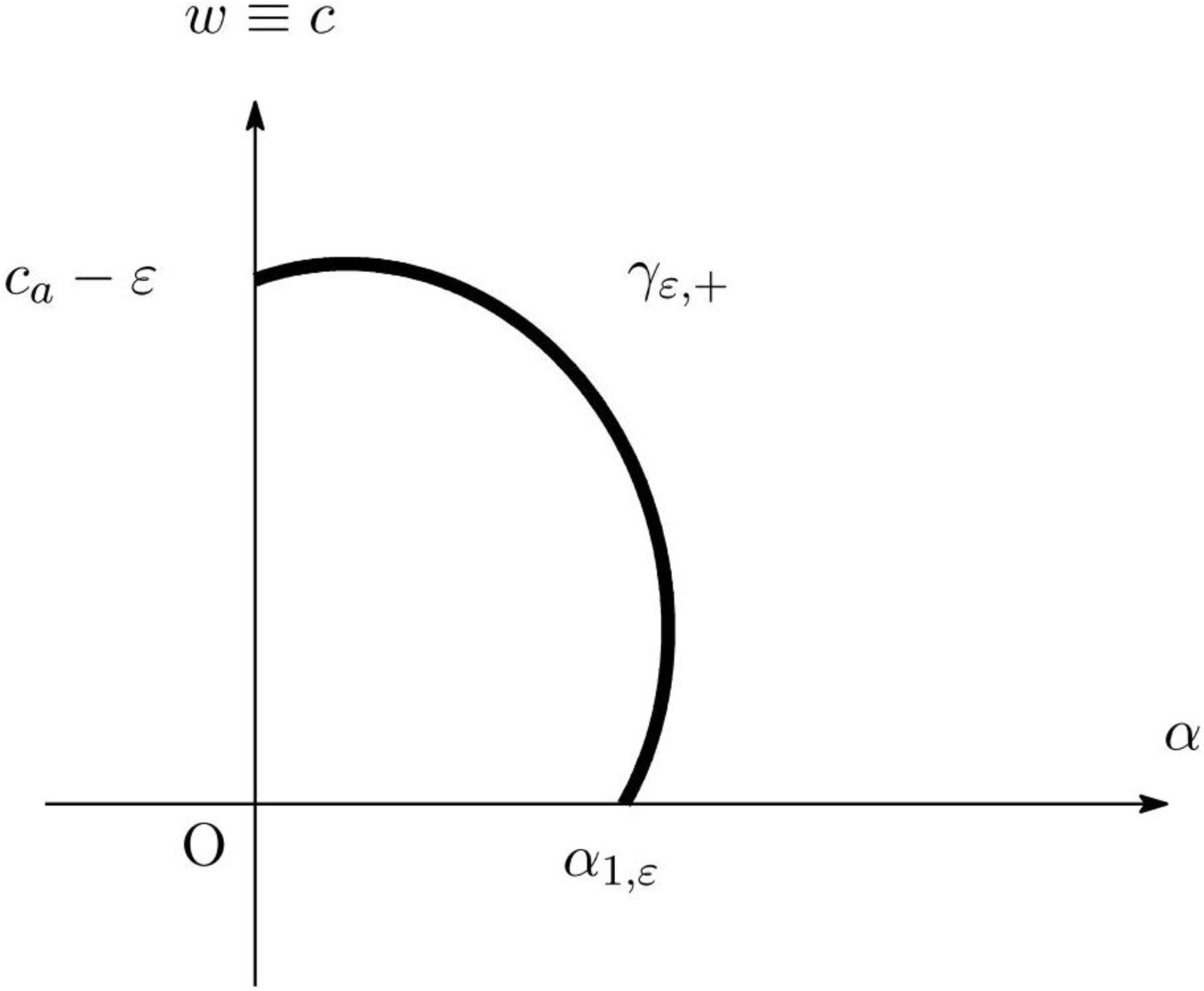}
} \centerline{(i) \hskip5.0cm (ii) }\caption{The component $\gamma
_{\varepsilon}$, and the bounded subcontinuum $\gamma_{\varepsilon, +}$.}%
\label{fig:double}%
\end{figure}

To prove the existence of $\gamma_{0}$ we shall employ Whyburn's topological
argument \cite[(9.12) Theorem]{Wh64}, applied to $\gamma_{0,\varepsilon}$,
$\gamma_{1,\varepsilon}$.
%Now, our analysis proceeds under the hypotheses of Propositions \ref{prop:bpc:e} and \ref{cor:bfie}, (A.0) and (A.2).
By Propositions \ref{cor:bfie} and \ref{prop:nbddalph}, we infer that
$\gamma_{0,\varepsilon}=\gamma_{1,\varepsilon}(=:\gamma_{\varepsilon})$ if
$\varepsilon> 0$ is small enough, see Figure \ref{fig:double}(i). More
precisely, Proposition \ref{prop:bpc:e} and \cite[Proposition 18.1]{Am76} tell
us that
\[
\gamma_{\varepsilon,+}:=\overline{\{(\alpha,w)\in\gamma_{\varepsilon}%
:\alpha>0\}},
\]
is a \textit{bounded} (compact) subcontinuum in $[0,\infty)\times
C^{1}(\overline{\Omega})$, satisfying
\begin{align}
\label{meet01}\gamma_{\varepsilon,+}\cap\Gamma_{0}=\{(\alpha_{1,\varepsilon
},\,0)\},\quad\gamma_{\varepsilon,+}\cap\Gamma_{1}=\{(0,\,c_{a}-\varepsilon
)\},
\end{align}
see Figure \ref{fig:double}(ii).
%
%\item $\mathcal{C}_\varepsilon \cap \Gamma_{0} = \{ (\lambda_{1, \varepsilon}, 0) \}$, i.e., $\mathcal{C}_\varepsilon$ does not meet $(\lambda, 0)$ for any $\lambda \neq 0, \lambda_{1,\varepsilon}$.
%\end{itemize}
%		r=-t^2+2t+3

%	% \begin{figure}[H]
%\begin{figure}[tbh]
%\begin{center}
%\includegraphics[scale=0.2]{fig18_0406gep02.eps}
%\end{center}
%\caption{The component $\gamma_{\varepsilon}$.}%
%\label{fig18_0406gep02}%
%\end{figure}

%	% \begin{figure}[H]
%\begin{figure}[tbh]
%\begin{center}
%\includegraphics[scale=0.2]{fig18_0406gep03.eps}
%\end{center}
%\caption{The bounded subcontinuum $\gamma_{\varepsilon, +}$.}%
%\label{fig18_0108a}%
%\end{figure}

%-----------------------

%\mpfns{{\bf added. the proof rewritten.}}
Now, let us analyze the limiting behavior of $\gamma_{\varepsilon,+}$ as
$\varepsilon\rightarrow0^{+}$.
%When $\varepsilon \to 0^{+}$, we see from Lemma \ref{prop:lam10} that $\alpha _{1,\varepsilon} \to0$ and from Proposition \ref{prop:bpc:e} that $c_a-\varepsilon\to c_a$. So,
We introduce the sets%
\begin{align*}
&  \liminf_{\varepsilon\rightarrow0^{+}}\gamma_{\varepsilon,+}:=\{(\alpha
,u)\in\mathbb{R}\times C^{1}(\overline{\Omega}):\lim_{\varepsilon \to 0^+
}\mathrm{dist}\,((\alpha,u),\gamma_{\varepsilon,+})=0\},\\
&  \limsup_{\varepsilon\rightarrow0^{+}}\gamma_{\varepsilon,+}:=\{(\alpha
,u)\in\mathbb{R}\times C^{1}(\overline{\Omega}):\liminf_{\varepsilon
\rightarrow0^{+}}\mathrm{dist}\,((\alpha,u),\gamma_{\varepsilon,+})=0\}.
\end{align*}
From the combination of \eqref{u:lbD} and \eqref{changeofv}, it follows that
given $\beta> 0$, there exists $\varepsilon_{0} > 0$ such that $w\geq
\beta^{\frac{1}{1-q}}u_{\beta,\varepsilon_{0}}$ on $\overline{D}$ for every
solution $w\in\mathcal{P}^{\circ}$ of $(R_{\alpha}^{\varepsilon})$ with
$\alpha\geq\beta$ and $0<\varepsilon\leq\varepsilon_{0}$. This implies that
$\alpha_{1,\varepsilon}<\beta$ if $0<\varepsilon\leq\varepsilon_{0}$, i.e.
$\alpha_{1,\varepsilon}\rightarrow0$ as $\varepsilon\rightarrow0^{+}$. In view
of \eqref{meet01},
\begin{equation}
(0,0),\ (0,c_{a})\in\liminf_{\varepsilon\rightarrow0^{+}}\gamma_{\varepsilon
,+}. \label{liminfga}%
\end{equation}
In a similar way as in \cite[Section 3]{RQUcpaa}, we can show that
$\bigcup\{\gamma_{\varepsilon,+}:\varepsilon>0\ \mbox{is small}\}$ is
precompact. Whyburn's topological argument \cite[(9.12) Theorem]{Wh64} can be
now applied to deduce that
\[
\gamma_{0,+}:=\limsup_{\varepsilon\rightarrow0^{+}}\gamma_{\varepsilon,+}%
\]
is a bounded subcontinuum in $[0,\infty)\times C^{1}(\overline{\Omega})$. In
addition, we infer
%\mpfns{modified}
from \eqref{liminfga} that $(0,0),\ (0,c_{a})\in\liminf_{\varepsilon
\rightarrow0^{+}}\gamma_{\varepsilon,+}\subseteq\gamma_{0,+}$.
%Now, by the same argument as in the proof of \cite[Theorem 1]{RQUcpaa}, we can prove that
%\begin{itemize}\setlength{\itemsep}{0.0cm}
%\item $\gamma_{0, +}$ consists of nonnegative solutions of $(R_{\alpha})$,
%\item  $\gamma_{0, \varepsilon}$ is nontrivial, \textit{i.e.,}  $\gamma_{0, \varepsilon} \not \subset \Gamma_{0}\cup\Gamma_1$,
%\item $\gamma_{0,+}$ satisfies \eqref{gcapGab}.
%\end{itemize}

Now, we verify that $\gamma_{0,+}$ consists of solutions of $(R_{\alpha})$.
Let $(\hat{\alpha},\hat{w})\in\gamma_{0,+}$. By definition, we can choose
$\varepsilon_{n}\rightarrow0^{+}$ and $(\alpha_{n},w_{n})\in\gamma
_{\varepsilon_{n}}$ such that $\alpha_{n}\rightarrow\hat{\alpha}\geq0$ and
$w_{n}\rightarrow\hat{w}$ in $C^{1}(\overline{\Omega})$. Since for all
$\varphi\in H^{1}(\Omega)$
%We shall investigate the limit $\hat{w}$ for the solution $w_{n}$ of $(R_{\alpha_n})$:%
\[
\int_{\Omega}\nabla w_{n}\nabla\varphi=\alpha_{n}\int_{\Omega}a(x)\left(
\frac{w_{n}}{w_{n}+\varepsilon_{n}}\right)  ^{1-q}w_{n}^{q}\varphi+\alpha
_{n}\int_{\partial\Omega}w_{n}\varphi,
\]
%When $\hat{w}(x)=0$, we see that
%\begin{align*}
%\left|  a(x) \left(  \frac{w_{n}(x)}{w_{n}(x) + \varepsilon_{n}} \right)  ^{1-q} w_{n}(x)^{q} \right|  \leq\| a \|_{C(\overline{\Omega})}|w_{n}(x)|^{q} \longrightarrow0.
%\end{align*}
%When $\hat{w}(x)>0$, we see that
%\begin{align*}
%a(x) \left(  \frac{w_{n}(x)}{w_{n}(x) + \varepsilon_{n}} \right)  ^{1-q} w_{n}(x)^{q} \longrightarrow a(x) \left(  \frac{\hat{w}(x)}{\hat{w}(x)+0}\right)^{1-q} \hat{w}(x)^{q} =  a(x)\hat{w}(x)^{q}.
%\end{align*}
%we use the Lebesgue convergence theorem to infer that
%\begin{align*}
%\int_{\Omega}a(x) \left(  \frac{w_{n}(x)}{w_{n}(x) + \varepsilon_{n}} \right)^{1-q} w_{n}(x)^{q} \varphi\longrightarrow\int_{\Omega}a(x) \hat{w}^{q} \varphi.
%\end{align*}
it follows by the Lebesgue dominated convergence theorem that
\begin{equation}
\int_{\Omega}\nabla\hat{w}\nabla\varphi=\hat{\alpha}\int_{\Omega}a(x)\hat
{w}^{q}\varphi+\hat{\alpha}\int_{\partial\Omega}\hat{w}\varphi,\quad
\forall\varphi\in H^{1}(\Omega). \label{hatuws}%
\end{equation}
Thus, $\hat{w}$ is a solution of $(R_{\hat{\alpha}})$ by elliptic regularity.

Next, we verify that $\gamma_{0,+}$ is nontrivial, \textit{i.e.,}
$\gamma_{0,+} \not \subset \Gamma_{0}\cup\Gamma_{1}$. Since $\gamma_{0,+}$ is
connected and joins $(0,0)$ to $(0, c_{a})$, the intermediate value theorem
shows that for $0< c< c_{a}$ we can pick $(\hat{\alpha}, \hat{w}) \in
\gamma_{0,+}$ such that $\hat{\alpha} \geq0$ and $\|\hat{w} \|_{C^{1}%
(\overline{\Omega})} = c$. We claim that $\hat{\alpha} > 0$, i.e.
$(\hat{\alpha}, \hat{w}) \not \in \Gamma_{0}\cup\Gamma_{1}$. To this end,
assume by contradiction that $\hat{\alpha}=0$. From the fact that
$(\hat{\alpha}, \hat{w}) \in\gamma_{0,+}$, we infer that there exist
$\varepsilon_{n} \to0^{+}$ and $(\alpha_{n}, w_{n}) \in\gamma_{\varepsilon
_{n}, +}$ such that $\alpha_{n} \to0^{+}$ and $w_{n} \to\hat{w}$ in
$C^{1}(\overline{\Omega})$. From \eqref{hatuws}, it follows that $\hat
{w}\equiv c$. However, from the definition of $w_{n}$ we obtain that
$\int_{\Omega}a (w_{n} + \varepsilon_{n})^{q-1}w_{n} + \int_{\partial\Omega}
w_{n} = 0$, and so, passing to the limit, that $c=c_{a}$, a contradiction.
%\[
%c^{q} \int_{\Omega}a + c \int_{\partial\Omega}b = 0 \quad(
%\mbox{implying $c = c_a$}).
%\]

Finally, we show how $\gamma_{0,+}$ meets $\Gamma_{0}$ and $\Gamma_{1}$. From
Proposition \ref{prop:bG00}(ii), we see that $\gamma_{0,+}$ does not meet any
point on $\Gamma_{1}$ except $(0,0)$ and $(0,c_{a})$. Moreover, Lemma
\ref{prop:glb} tells us that $\gamma_{0,+}$ does not meet $\Gamma_{0}$, so
that $\gamma_{0,+}$ satisfies \eqref{gcapGab}.

To sum up, $\gamma_{0}:=\gamma_{0,+}$ is as desired.
%
%assuming the conditions in Propositions \ref{prop:bfzero}, \textit{i.e.,} assuming $(A.1)$ and $(A.3)$, we deduce that $\gamma_{0,+}$ does not meet $\Gamma_{0}$, which implies that $\gamma_{0,+}$ satisfies \eqref{gcapGab}. To sum up, $\gamma_{0}:=\gamma_{0,+}$ is as desired. Indeed, that if $(\alpha, u)  \in\gamma_{0,+}\setminus\{ (0,0), (0, c_a) \}$, then we have that $\alpha> 0$, and $u$ is a \textit{nontrivial}  nonnegative solution of $(R_{\alpha})$.
%---------------------------
%The rest of the argument is quite similar to that in the proof of \cite[Theorem 1]{RQUcpaa}, and so, we omit it. But, it should be remarked that our current situation is just the same as case (a) illustrated in \cite[Figures 4 and 5]{RQUcpaa}.
Indeed, assertion (i) follows from Proposition \ref{cor:Aal:ordered}. The
exactness assertion in (ii) comes from Theorem \ref{prop:exactm}. The
positivity assertion in (iii) is a consequence of assertion (i), the
boundedness assertion follows from Proposition \ref{cor:bounds} and
Proposition \ref{alpha}, and finally, \eqref{hatgam0} follows from the second
assertion of Lemma \ref{prop:glb}. The proof is now complete. \qed\newline

\begin{rem}
\label{4.5}Assuming only $(A.0)$ and $(A.3)$ we can establish, for
\textit{any} $q\in(0,1)$, the existence of a subcontinuum $\gamma
_{0}=\{(\alpha,w)\}$ in $[0,\infty)\times C^{1}(\overline{\Omega})$ of
solutions of $(R_{\alpha})$ satisfying \eqref{gcapGab} and $w>0$ in
$\Omega_{+}^{a}\cup D\cup\Gamma_{\partial\Omega}$ whenever $(\alpha
,w)\in\gamma_{0}$. Indeed, if $(\alpha,w)\in\gamma_{0}$ for some $\alpha>0$
then there exist $\varepsilon_{n}\rightarrow0^{+}$, $\alpha_{n}\rightarrow
\alpha$, and $w_{n}\rightarrow w$ in $C^{1}(\overline{\Omega})$ such that
$(\alpha_{n},w_{n})\in\gamma_{\varepsilon_{n},+}$, implying that $w_{n}%
\in\mathcal{P}^{\circ}$ is a solution of $(R_{\alpha_{n}}^{\varepsilon_{n}})$.
Applying a sub and supersolutions argument as in the proof of Proposition
\ref{prop:nbddalph} with $u_{n}:=\alpha_{n}^{1/(q-1)}w_{n}$, we deduce that
$w>0$ in $\Omega_{+}^{a}\cup D\cup\Gamma_{\partial\Omega}$. Note that
Proposition 3.6 still holds for solutions of $(R_{\alpha})$ that are positive
in $\Omega_{+}^{a}$, so the component $\hat{\gamma}_{0}=\{(\alpha,w)\}$ of
solutions of $(R_{\alpha})$ that includes $\gamma_{0}$ has the same nature as
in Theorem \ref{maint}(iii), but is composed now by solutions that are
positive in $\Omega_{+}^{a}\cup D\cup\Gamma_{\partial\Omega}$. \newline
\end{rem}

\noindent\textit{Proof of Theorem \ref{thm:main}.}

First we verify (i). The assertion $\alpha_{s}\in\left(  0,\infty\right)  $
follows from Theorem \ref{thm:curve}(ii) and Propositions
\ref{cor:Aal:ordered} and \ref{alpha}, whereas the second assertion follows
from Proposition \ref{prop:exnon}, thanks to Theorem \ref{tpos}. Assertion
(ii) is deduced from Theorem \ref{thm:curve}, Proposition \ref{prop:bG00} and
Theorem \ref{prop:exactm}. Indeed, $\mathcal{C}_{1}$ is given by Proposition
\ref{prop:bG00}. Finally, the existence and properties of the component
$\mathcal{C}_{\ast}$ in (iii) are proved by combining Theorem \ref{thm:curve}%
(ii) and Theorem \ref{maint}. \qed\newline

\end{document}